\documentclass{article}

\usepackage{arxiv}

\usepackage{amssymb}
\usepackage{amsthm}
\usepackage{amsmath}
\usepackage{color,soul}
\usepackage[dvipsnames]{xcolor}
\usepackage{graphicx}
\usepackage{tcolorbox}
\usepackage{ulem}
\usepackage{appendix}

\usepackage{mathtools}

\newtheorem{theorem}{Theorem}
\newtheorem{lemma}{Lemma}
\newtheorem{corollary}{Corollary}

\theoremstyle{definition}
\newtheorem{example}{Example}
\newtheorem{remark}{Remark}

\allowdisplaybreaks[4]

\begin{document}
	
\title{Kullback-Leibler Divergence and Akaike Information Criterion in General Hidden Markov Models}

\author{
  Cheng-Der Fuh
  \\
  Department of Statistics\\
  Zhejiang University City College\\
  Hangzhou 310015, China \\
  \texttt{cdffuh@gmail.com} \\
   \And
 Chu-Lan Michael Kao \\
  Institute of Statistics\\
  National Yang Ming Chiao Tung University\\
  Hsinchu 30010, Taiwan \\
  \texttt{chulankao@gmail.com} \\
   \And
    Tianxiao Pang\\
    School of Mathematical Sciences\\
    Zhejiang University\\
    Hangzhou 310058, China\\
    \texttt{txpang@zju.edu.cn}
}

\maketitle

\begin{abstract}
To characterize the Kullback-Leibler divergence and Fisher information in general parametrized hidden Markov models, in this paper, we first show that the log likelihood and its derivatives can be represented as an additive functional of a Markovian iterated function system, and then provide explicit characterizations of these two quantities through this representation. Moreover, we show that Kullback-Leibler divergence can be locally approximated
by a quadratic function determined by the Fisher information. Results relating to the Cram\'{e}r-Rao lower bound and the H\'{a}jek-Le Cam local asymptotic minimax theorem are also given. As an application of our results, we provide a theoretical justification of using Akaike information criterion (AIC) model selection in general hidden Markov models. Last, we study three concrete models: a Gaussian vector autoregressive-moving average  model of order $(p,q)$, recurrent neural networks, and
temporal restricted Boltzmann machine, to illustrate our theory.
\end{abstract}

\keywords{
AIC, Boltzmann machine, Cram\'{e}r-Rao lower bound, Fisher information, H\'{a}jek-Le Cam theorem, hidden Markov model, Kullback-Leibler divergence, Markovian iterated function system, recurrent neural network.}

\section{Introduction}\label{sec:intro}

Kullback-Leibler (KL) divergence, also called relative entropy, has been widely used in information theory, machine learning, statistics, econometrics, and others. Its applications include
information theory (\cite{aghamohammadi2018extreme}), speech recognition via deep neural networks (\cite{yu2013kl}), chemical kinetics (\cite{hangos2010engineering}),
physics (\cite{beck2009generalised}), statistics and econometrics (\cite{White_1982, Gourieroux_Monfort_Trognon_1984, kitamura1997information,
	Imbens_Spady_Johnson_1998}).
Theoretical properties of the KL-divergence and its relationship with the Fisher information matrix have been well established, particularly for models with independent and identically distributed (i.i.d.) observations.

Nonetheless, many real applications now build on more complex hidden Markov models (HMMs) with finite states, or even general hidden Markov models (GHMMs) with general states. The former includes machine learning applications in speech recognition (\cite{juang1991hidden}) and computational biology (\cite{marioni2006biohmm}), econometric applications with Markov switching models (\cite{hamilton1989new, Calvet_Fisher_2001}), Markov switching GARCH models (\cite{Cai_1994, Hamilton_Susmel_1994}) and many other applications. The latter includes factorial HMMs (\cite{ghahramani1997factorial}), switching state-space models (\cite{Kim_1994, Kim_Nelson_1998, ghahramani2000variational}) and adversarial models (\cite{goodfellow2014generative}) in machine learning,  (G)ARCH models (\cite{engle1982autoregressive, bollerslev1986generalized, fan2003nonlinear, hall2003inference, francq2012strict}) and stochastic volatility models (SVs) (\cite{clark1973subordinated, taylor1986modeling, barndorff2004power}) in statistics and econometrics, and others from various disciplines. It is also known that KL-divergence plays an important role in HMMs and GHMMs. For example, \cite{smith2006markov} applies KL-divergence for model selection in Markov switching models. \cite{fuh2003} and \cite{Andreoua_Ghysels_2008} use KL-divergence to detect change points for HMMs, while \cite{fuh2021} studies the detection for GHMMs. \cite{FuhMei2015} further provides a numerical computational method via Fredholm integral equations in a two-state HMM. See also \cite{gorban2010entropy} and \cite{travers2014exponential}, as well as \cite{obremski2020complexity} for the more general R\'{e}nyi entropy in Markov models.
This motivates us to have a theoretical investigation of the KL-divergence in GHMMs.

Note that there are many results and mathematical mechanisms for i.i.d.\ models, but many of them cannot be directly applied to GHMM, or even to HMM. The main reason is that the log likelihood of a HMM or GHMM is {\it not} the sum of i.i.d.\ random variables or even a functional of Markov chains; thus the classical law of large numbers (LLN) approach cannot be directly applied. Instead, \cite{Leroux:1992} applies Kingman's subadditive ergodic theorem to provide a generalized KL-divergence, while \cite{bickel1998asymptotic} uses an ergodic process to approximate the log likelihood function in a finite-state HMM. \cite{douc2011consistency} applies the Shannon-Breiman-McMillan theorem  to have the limit as the KL-divergence in a general-state HMM. However, their results require stationarity for the HMM, and do not characterize the KL-divergence in HMM or GHMM, nor does the Fisher information; hence many asymptotic properties for HMM and GHMM remain uninvestigated, including KL-divergence and its relationship with Fisher information.

To formally explain this phenomenon in details, we first follow the definition in \cite{Fuh2006} to define the GHMM. Let $\{X_n, n \geq 0 \}$ be a Markov
chain on a general state space $\mathcal{X}$, with transition probability kernel $p_{\theta}(x,\cdot)= P^{\theta}\{X_1 \in \cdot|X_0=x\}$ and stationary
probability $\pi(\cdot):=\pi_{\theta}(\cdot)$ with respect to a $\sigma$-finite measure~$Q$ on $\mathcal{X}$, where $\theta \in \Theta \subseteq {\bf R}^q $
denotes the unknown parameter. Let $Y_{0:n}$ be the observations from $Y_0$ to $Y_n$ such that $Y_n\in {\bf R}^d$ with a distribution depending on $X_n$ and $Y_{n-1}$, but independent to others. Let $f(\cdot; \theta|x, y)$ be the probability density function (pdf) of $Y_n$ given $X_n=x$ and $Y_{n-1}=y$, with respect to a $\sigma$-finite measure~$\tilde{Q}$ on ${\bf R}^d$. Further let $f(\cdot; \theta|x_0)$ be the pdf of $Y_0$ given $X_0 = x_0$. Note that this setting includes
interesting examples such as Markov-switching autoregression models, (G)ARCH models, stochastic volatility models, recurrent neural networks (RNNs) and temporal restricted Boltzmann machine. When $\mathcal{X}$ is a finite state space and $Y_n$ are independent for given $X_{n}$, this is the classical hidden Markov model.

For given random observations $Y_{0:n}$, the full likelihood  is
\begin{align}\label{eqn:lik}
L(\theta; Y_{0:n})
= & \int_{x_0 \in \mathcal{ X}} \cdots  \int_{x_n \in \mathcal{ X}} \pi_\theta(x_0)
f(Y_0;\theta|x_0) \\
\notag
& \times \prod_{t=1}^n p_\theta(x_{t-1},x_t) f(Y_t;\theta| x_t, Y_{t-1})
Q(dx_n)\cdots Q(dx_0).
\end{align}
In addition, denote $\ell(\theta; Y_{0:n}) := \log L(\theta; Y_{0:n})$ as the $\log$
likelihood.  Note that in (\ref{eqn:lik}) the initial distribution of $X_0$ is taken as the stationary distribution $\pi_\theta(\cdot)$ for convenience, indeed any suitable initial distribution $\bar{\nu}(\cdot)$ works well.

Then, for any two parameters $\theta_0$~and~$\theta_1$, the KL-divergence $K(\theta_1, \theta_0)$ is defined as
\begin{align}\label{KL-HMM}
    K(\theta_1, \theta_0) =  \lim_{n \rightarrow \infty}\frac{1}{n}\left[ \ell(\theta_1; Y_{0:n}) - \ell(\theta_0; Y_{0:n})\right],~~~P^{\theta_1}\mbox{-a.s.},
\end{align}
where $P^{\theta}$ denotes the probability measure when $(Y_0, \cdots, Y_n)$
are distributed according to $L(\theta; \cdot)$. In addition, the Fisher
information under $P^{\theta_0}$ can be defined as
\begin{equation}\label{Fisher-Bickel}
I(\theta_0) = -\lim_{n \rightarrow \infty} \frac{1}{n} \frac{\partial^2 \ell(\theta_0; Y_{0:n})}{\partial \theta_0\partial \theta_0^t},~~~P^{\theta_0}\mbox{-a.s.},
\end{equation}
where the superscript $t$ denotes the transpose; see the last equation on page~2047 of \cite{Fuh2006}.

When $\{Y_n, n \geq 0\}$ are i.i.d. random variables  with pdf $f(y;\theta|x,y_0) = f(y;\theta)$, then
\begin{equation}\label{eqn:logLik-iid}
    \ell(\theta; Y_{0:n}) = \sum_{t=0}^n \log f(Y_t; \theta)
\end{equation}
is a sum of i.i.d.\ random variables $\{\log f(Y_t;\theta), t \geq 0\}$.
Hence, under some regularity conditions, by \eqref{KL-HMM} and the strong law of large numbers (SLLN), we have
\begin{align}\label{KL-iid}
    K(\theta_1, \theta_0) = E^{\theta_1}\left[ \log f(Y_1;\theta_1) - \log f(Y_1;\theta_0) \right] = \int_{{\bf R}^d} \log \frac{f(y;\theta_1)}{f(y;\theta_0)}f(y;\theta_1) \tilde{Q}(dy),
\end{align}
where $E^{\theta}$ denotes the expectation under $P^\theta$.
Similarly, for $i,j=1,\cdots,q$,
\begin{equation}\label{eqn:2ndDeriv-iid}
   \frac{\partial^2 \ell(\theta; Y_{0:n})}{\partial \theta_i \partial \theta_j} = \sum_{t=0}^n \frac{\partial^2 \log f(Y_t; \theta)}{\partial \theta_i \partial \theta_j}
\end{equation}
is a sum of i.i.d.\ random variables $\{ \frac{\partial^2 \log f(Y_t; \theta)}{\partial \theta_i \partial \theta_j}, t \geq 0\}$, therefore by \eqref{Fisher-Bickel} and SLLN, we have
\begin{equation}\label{Fisher-iid}
[I(\theta)]_{i,j} = -E^{\theta}\left[ \frac{\partial^2 \log f(Y_1; \theta)}{\partial \theta_i \partial \theta_j}\right].
\end{equation}
Finally, with the help of \eqref{KL-iid} and \eqref{Fisher-iid}, it is known
that as $\theta_1 \rightarrow \theta_0$, we have
\begin{equation}\label{KL-converge}
    K(\theta_1, \theta_0) = (\theta_1 - \theta_0)^t \frac{I(\theta_0)}{2} (\theta_1-\theta_0) + O(\Vert \theta_1 - \theta_0 \Vert^3),
\end{equation}
where $\Vert \cdot \Vert$ denotes the Euclidean norm.

Nevertheless, for HMM and GHMM cases, we do not have
\eqref{eqn:logLik-iid}, which precludes us from directly obtaining \eqref{KL-iid}
through SLLN. Similarly, since we do not have \eqref{eqn:2ndDeriv-iid},
\eqref{Fisher-iid} cannot be derived using the same argument. As a consequence,
although \eqref{KL-converge} has been long conjectured in the literature (see,
for example, Remark~2 in \cite{fuh2004bahadur}), a rigorous proof is still lacking.

Note that these difficulties are all highly related to the complex
structure of \eqref{eqn:lik}.
Therefore, in this paper, we use an innovative
representation of the log likelihood $\ell(\theta; Y_{0:n})$
and its derivatives in GHMM, which
gets around this complexity. By such, we provide
characterizations of the KL-divergence and Fisher information, and prove
the relationship between these two via the corresponding convergence in
\eqref{KL-converge}.

Given these newly developed characterizations, we further provide the
Cram\'{e}r-Rao lower bound and H\'{a}jek-Le Cam local asymptotic minimax
theorem (\cite{van2002statistical}) for GHMM, which shows that the classical
bounds in i.i.d.\ scenarios remain valid for GHMM. We also show that as in the
i.i.d.\ case, the KL-divergence satisfies the non-negativity and additivity
properties. However, it is not convex in general, which is in contrast to the
traditional i.i.d.\ or Markov chain cases for which the KL-divergence is convex.
As another application of our results, we further provide a theoretical justification of using
Akaike information criterion (AIC) model selection in GHMMs.

The rest of the paper is organized as follows. In Section~\ref{sec:main} we present conditions and state main results. Section \ref{sec:AIC} studies the application to AIC model selection. To illustrate our theoretical results, three concrete models: a Gaussian vector autoregressive-moving average  model of order $(p,q)$, recurrent neural networks, and
temporal restricted Boltzmann machine,  are discussed in Section \ref{sec:Examples}.
Section~\ref{sec:conclude} concludes.
All proofs of theoretical results are given in Appendix.

\section{Main Results}\label{sec:main}

We split this section into three parts. Section \ref{subsec:condi} defines notations and states conditions. Section \ref{subsec:KeyLemma} presents preliminary results, which show that the log likelihood and the
derivatives of the log likelihood can be represented as an additive functional of a Markovian iterated function system (MIFS). Section \ref{subsec:MainResults} states our main results, which include characterizations of the KL-divergence, Fisher information matrix, and the relationship between
these two for a GHMM. Moreover, we show the results relating to the Cram\'{e}r-Rao lower bound and  H\'{a}jek-Le Cam local asymptotic minimax theorem.

\subsection{Notations and Conditions}\label{subsec:condi}

Denote $E_x^{\theta}$ as the expectation defined under
$P^{\theta}$ with initial state $X_0=x,$ and $E_{(x,y)}^\theta$ as
the expectation defined under $P^{\theta}$ with initial state
$(X_0,Y_0)=(x,y)$. For any $1 \leq i \leq q$ and positive integer~$k$, let
$D_i$ be the partial derivative with respect to the $i$-th dimension of
$\theta$ in some neighborhood $N_\delta(\theta_0) := \{\theta: \Vert \theta -
\theta_0\Vert < \delta\}$ of the true value $\theta_0$, and let $(D_i)^k$ be the
corresponding $k$-th partial derivative. In addition, for a given non-negative
integer vector $\nu=(\nu^{(1)},\cdots,\nu^{(q)})$, write
$|\nu|=\nu^{(1)}+\cdots+\nu^{(q)}$, $\nu!=\nu^{(1)}! \cdots \nu^{(q)}!$,
and let $D^{\nu}_\theta:=D^{\nu}=(D_1)^{\nu^{(1)}} \cdots (D_q)^{\nu^{(q)}}$ denote the
$\nu$-th derivative with respect to $\theta$ in $N_\delta(\theta_0)$.

The following conditions will be used throughout the rest of this paper.

C1. For a given $\theta \in \Theta$, the Markov chain $\{(X_n, Y_n), n \geq 0 \}$
is aperiodic, irreducible, and satisfies
\begin{align*}
\lim_{n \rightarrow \infty} \sup_{\substack{x \in \mathcal{X}, y \in {\bf R}^d, |h| \leq w}} \left\vert \frac{E_{(x,y)}^\theta[h(X_n, Y_n)] - \int
  h(s)\pi(ds)}{w(x,y)} \right\vert =0,
\end{align*}
\begin{align*}
\sup_{(x,y) \in \mathcal{X} \times {\bf R}^d} & \frac{E_{(x,y)}^\theta[w(X_p, Y_p)]}{w(x,y)}  < \infty,
\end{align*}
with some weight function $w(\cdot, \cdot)$ and $p \geq 1$.
Assume that
\begin{equation}\label{C1-3}
0< p_{\theta}(x_0,x_1)< \infty ~~ \mbox{for all} ~~ x_0,x_1 \in \mathcal{ X},
\end{equation}
and
\begin{equation}\label{C1-4}
0< \sup_{x \in \mathcal{ X}}
f(y_1;\theta|x,y_0)< \infty ~~ \mbox{for all} ~~ y_0,y_1 \in {\bf R}^d.
\end{equation}
Since $Q$ is $\sigma$-finite, there exist pairwise disjoint $\mathcal{ X}_n$'s such that  $\mathcal{ X}= \cup_{n=1}^{\infty} \mathcal{ X}_n$,
and $0 < Q(\mathcal{ X}_n) < \infty$. Assume that
\begin{equation}\label{C1-5}
E^\theta \left[\sum_{n=1}^{\infty} \frac{1}{2^n}\sup_{x \in \mathcal{ X}_n} f(Y_1;\theta|x,y_0)\right]
< \infty ~~ \mbox{ for all } y_0 \in {\bf R}^d.
\end{equation}
Furthermore, let
\begin{equation*}
\tilde{f}_{\theta}(y_0,y_1)= \sup_{x_0 \in \mathcal{ X}} \int_{x \in \mathcal{ X}} p_{\theta}(x_0,x)
f(y_1;\theta|x,y_{0}) Q(dx),
\end{equation*}
and assume that there exists $p \geq 1$ such that
\begin{align}
& \sup_{(x_0, y_0) \in \mathcal{ X} \times {\bf R}^d} E^{\theta}_{(x_0,y_0)} \left\{ \log \bigg( (\tilde{f}_{\theta}(y_0,Y_1))^p
\frac{w(X_p,Y_p)}{w(x_0,y_0)}  \right) \bigg\}  < 0 \label{C1-6}, \\
& \sup_{(x_0,y_0) \in \mathcal{ X} \times {\bf R}^d} E^{\theta}_{(x_0,y_0)} \left\{ \tilde{f}_{\theta}(y_0,Y_1) \frac{w(X_1,Y_1)}{w(x_0,y_0)} \right\}  < \infty. \label{C1-7}
\end{align}

C2. The true parameter~$\theta_0$ is an interior point of $\Theta$.
For all $x\in \mathcal{ X}$, $y_0,y_1 \in {\bf R}^d$, $\theta \in \Theta
\subset {\bf R}^q$ and $\nu$ with $|\nu| \leq r$, the partial derivatives
$D^\nu f(y_0;\theta|x)$ and $D^\nu f(y_1;\theta|x,y_0)$ exist. In addition,
for all $x_0,x \in \mathcal{ X}$, $\theta \mapsto p_\theta(x_0,x)$ and
$\theta \mapsto \pi_{\theta}(x_0)$ have $r$th-order continuous derivatives in
some neighborhood $N_\delta(\theta_{0})$ of $\theta_0$.

C3. For all $\nu$ with $|\nu| \leq r$ and $x_0 \in \mathcal{X}$
\begin{equation*}
    \int_{x \in \mathcal{X}} \sup_{\theta \in N_\delta(\theta_0)} \left\vert
D^\nu \pi_\theta(x)\right\vert Q(dx) < \infty
\end{equation*}
and
\begin{equation*}
\int_{x \in \mathcal{X}} \sup_{\theta \in N_\delta(\theta_0)} \left\vert D^\nu p_{\theta}(x_0,x)\right\vert Q(dx) < \infty.
\end{equation*}

C4. For all $x \in \mathcal{ X}$, $y_0 \in {\bf R}^d$ and $\theta \in \Theta$,
\begin{align*}
E_x^{\theta}| D^{\nu} f(Y_0;\theta|x)|^r < \infty, ~~
E_{(x,y_0)}^{\theta}| D^{\nu}f(Y_1;\theta|x,y_0)|^r < \infty
\end{align*}
for $1 \leq |\nu| \leq r$, and
\begin{align*}
E_x^{\theta} \left( \sup_{\theta \in N_\delta(\theta_0)} | D^{\nu} f(Y_0;\theta|x)|^r \right) < \infty, \\
E_{(x,y_0)}^{\theta} \left( \sup_{\theta \in N_\delta(\theta_0)} | D^{\nu}f(Y_1;\theta|x,y_0)|^r \right) < \infty
\end{align*}
for $|\nu| = r+1$.

C5.
\begin{equation*}
   E^{\theta_0}\left( \sup_{\Vert\theta-\theta_0\Vert<\delta} \sup_{x_0,x_0',x_1, x_1' \in \mathcal{X}} \frac
  {f(Y_0; \theta|x_0) f(Y_1; \theta|x_1, Y_0)}
  {f(Y_0; \theta|x_0') f(Y_1; \theta|x_1', Y_0)}\right)^r < \infty.
\end{equation*}

C6. For any $\theta \in N_\delta(\theta_0)$ and $\nu$ with $|\nu| \leq r$,
\begin{equation*}\label{C6-3}
\left\vert D^\nu p_{\theta}(x_0,x_1) \right\vert < \infty ~~ \mbox{for all} ~~ x_0,x_1 \in \mathcal{ X},
\end{equation*}
\begin{equation*}\label{C6-4}
\sup_{x \in \mathcal{ X}}\left\vert D^\nu
f(y_1;\theta|x,y_0) \right\vert < \infty ~~ \mbox{for all} ~~ y_0,y_1 \in {\bf R}^d,
\end{equation*}
\begin{equation*}\label{C6-5}
E^\theta \left[\sum_{n=1}^{\infty} \frac{1}{2^n}\sup_{x \in \mathcal{ X}_n} \left\vert D^\nu f(Y_1;\theta|x,y_0) \right\vert \right]
< \infty ~~ \mbox{ for all } y_0 \in {\bf R}^d.
\end{equation*}
Furthermore, let
\begin{equation*}
\tilde{f}_{\theta}^\nu(y_0,y_1)= \sup_{x_0 \in \mathcal{ X}} \int_{x \in \mathcal{ X}} D^\nu \left\{ p_{\theta}(x_0,x)
f(y_1;\theta|x,y_{0}) \right\} Q(dx),
\end{equation*}
and assume that there exists $p \geq 1$ such that
\begin{align*}
& \sup_{(x_0, y_0) \in \mathcal{ X} \times {\bf R}^d} E^{\theta}_{(x_0,y_0)} \left\{ \log \left( \left\vert \tilde{f}_{\theta}^\nu(y_0,Y_1) \right\vert^p
\frac{w(X_p,Y_p)}{w(x_0,y_0)}  \right) \right\}  < 0, \\
& \sup_{(x_0,y_0) \in \mathcal{ X} \times {\bf R}^d} E^{\theta}_{(x_0,y_0)} \left\{ \left\vert \tilde{f}_{\theta}^\nu(y_0,Y_1) \right\vert
\frac{w(X_1,Y_1)}{w(x_0,y_0)} \right\}  < \infty.
\end{align*}

\begin{remark}
Conditions C1 and C2--C5 are essentially the same as conditions C1 and
C2'--C5'
in \cite{Fuh2006}, respectively. The purpose of the additional condition C6, on the
other hand, is to extend \eqref{C1-3}--\eqref{C1-7} in C1 to higher-order
derivatives in some neighborhood of~$\theta_0$. Many commonly used models satisfy these conditions, including Markov switching models, ARMA models, (G)ARCH models as well as stochastic volatility models; see  \cite{Fuh2006} for details. Furthermore, we will check conditions C1--C6 also hold under RNN and temporal restricted Boltzmann machine with specific distributions.

\end{remark}

\subsection{Preliminary Results}\label{subsec:KeyLemma}

\cite{Fuh2006} has represented the log likelihood $\ell(\theta; \cdot)$ as an additive functional of a MIFS as follows. To be more specific, we consider
the function space
\begin{align*}
    {\bf M} = \bigg\{  h \Big\vert h:\mathcal{ X} \mapsto {\bf R} ~\textrm{is~}Q
    \textrm{-measurable},  \int_{x \in \mathcal{X}} |h(x)| Q(dx) < \infty
      \textrm{~and~}\sup_{x \in \mathcal{X}} |h(x)| < \infty \bigg\}.
    \end{align*}
Moreover, for $t=1,\cdots,n$, define the random functions ${\bf
P}_\theta(Y_0)$ and ${\bf P}_\theta(Y_j)$ on $(\mathcal{ X} \times {\bf R}^d)
\times {\bf M}$ as
\begin{align}
\notag
& {\bf P}_\theta(Y_0)h(x) = \int_{x_0 \in \mathcal{ X}} f(Y_0;\theta|x_0) h(x_0)
Q(dx_0), ~~~\textrm{a~constant~functional,}
\\
\notag
& {\bf P}_\theta(Y_t)h(x) = \int_{s \in \mathcal{ X}} p_{\theta}(s,x)
f(Y_t;\theta|x, Y_{t-1}) h(s) Q(ds),
\end{align}
and define the composition of two random functions as
\begin{align*}
& {\bf P}_\theta(Y_{t+1})\circ{\bf P}_\theta(Y_{t})h(x) \\
= &\int_{z \in \mathcal{ X}} p_\theta(z,x) f(Y_{t+1};\theta|x, Y_{t})\times \bigg(\int_{s \in \mathcal{ X}} p_\theta(s,z) f(Y_{t};\theta|z,Y_{t-1}) h(s) Q(ds)\bigg) Q(dz).
\end{align*}
Now, consider
\begin{eqnarray}\label{mn}
M_n:= {\bf P}_\theta(Y_n) \circ \cdots \circ {\bf P}_\theta(Y_1) \circ {\bf P}_\theta(Y_0).
\end{eqnarray}
Further denote $ \langle h \rangle := \int_{x \in \mathcal{ X}} h(x) Q(dx)$.
Then, we have
\begin{align}\label{function_g}
\ell(\theta, Y_{0:n}) & = \log L(\theta;Y_{0:n}) = \log \langle M_n \pi \rangle  \notag\\
& = \sum_{t=1}^n \log \frac{\langle M_t \pi \rangle}{\langle M_{t-1} \pi \rangle} + \log \langle M_0 \pi \rangle \notag\\
& =: \sum_{t=1}^n g^0( M_t^0, M_{t-1}^0) +g_0^0(M_0^0),
\end{align}
where
\begin{align}\label{function_gg}
g^0( M_t^0, M_{t-1}^0)=\log \frac{\langle M_t \pi \rangle}{\langle M_{t-1} \pi \rangle},\quad g_0^0(M_0^0)=\log \langle M_0 \pi \rangle.
\end{align}
In other words, $\ell(\theta)$ is an additive functional of  $\{((X_n, Y_n), M_n), n \geq 0\}$.
In addition, \cite{Fuh2006} shows that $\{((X_n, Y_n), M_n), n \geq 0\}$ forms
an ergodic Markov chain, induced by the MIFS based on \eqref{mn}, on the state
space $(\mathcal{X} \times {\bf R}^d) \times {\bf M}$. \cite{Fuh2006} further uses this
result to prove the SLLN for the log likelihood. The rate of convergence of $\{((X_n, Y_n), M_n), n \geq 0\}$ to its invariant measure is studied in \cite{fuh2021a}.

The following lemmas extend this idea to the derivatives of $\ell(\theta;\cdot)$. To do so, for any $q$-dimensional non-negative integer vector $\nu = (\nu^{(1)}, \cdots, \nu^{(q)})$, define
\begin{eqnarray*}
    W_n^\nu = D^\nu M_n =
     (D_1)^{\nu^{(1)}} \cdots (D_q)^{\nu^{(q)}} (M_n).
\end{eqnarray*}
Now let us consider all derivatives with order~$r$ or less. Note that for a
fixed integer $r \geq 1$,  there are exactly $K = {(r+q)!}/(r!q!)$ different
$\nu$ satisfying $|\nu| \leq r$. Label all such $\nu$ by $\nu_1, \nu_2, \cdots,
\nu_K$, and let $W_n^{(r)} = (W_n^{\nu_1}, W_n^{\nu_2}, \cdots,
W_n^{\nu_K})^t$.

The first lemma shows that we can construct a MIFS through $W_n^{(r)}$.
\begin{lemma}\label{thm:rep}
Assume conditions C1--C6 hold with some $r \geq 1$. Then, for any $\theta \in N_\delta(\theta_0)$,
\begin{equation*}
\{((X_n, Y_n), W_n^{(r)}), n \geq 0\}
\end{equation*}
is an aperiodic, $(\mathcal{ X} \times {\bf R}^d)
\times {\bf M}^K$-irreducible and Harris-recurrent Markov chain.
\end{lemma}
\noindent See the supplementary for the proof.

The second lemma shows that the derivatives of $\ell(\theta; \cdot)$ can be represented as an additive functional of this particular MIFS.
\begin{lemma}\label{thm:additiv}
Assume conditions C1--C6 hold with some $r \geq 1$. Then, for any $\theta \in N_\delta(\theta_0)$ and any
$q$-dimensional non-negative integer vector $\nu$ with $|\nu| \leq r$, there exists function $g^\nu$ and $g_0^\nu$ such that
\begin{eqnarray}\label{DL-to-g}
	D^\nu \ell(\theta;Y_{0:n}) = \sum_{t=1}^n g^\nu (W_{t}^{(|\nu|)}, W_{t-1}^{(|\nu|)}) + g_0^\nu(W_0^{(|\nu|)}).
\end{eqnarray}
\end{lemma}
\noindent See the supplementary for the proof.

Lemmas \ref{thm:rep} and \ref{thm:additiv} are almost the same as Lemmas 3 and 5 in \cite{fuh_pang_2022},
respectively, for a two-layer HMM, we include them here for completeness.
Combining Lemmas \ref{thm:rep} and \ref{thm:additiv}, we can apply the LLN for MIFS to evaluate $D^\nu \ell(\theta;\cdot)$. This further leads to the main results in the next subsection.

\subsection{Main Results}\label{subsec:MainResults}

We will use Lemmas \ref{thm:rep} and \ref{thm:additiv} to evaluate Fisher information, KL-divergence and other properties. However, as these quantities might involve different probability measures as well as $\ell(\theta;\cdot)$ evaluated at different $\theta$, some additional notations are needed to clarify the statement. For $i=0,1$, let $W_{n,\theta_i}^{(r)}$ be the $W_n^{(r)}$ constructed with the $\ell(\theta; Y_{0:n})$ evaluated at $\theta=\theta_i$. In addition, for any $1 \leq j,k \leq q$, let $I_{jk}(\theta_0)$ be the $(j,k)$-th component
in the Fisher information matrix $I(\theta_0)$. Further denote $\vec{0} =
(0,0,\cdots, 0) \in {\bf R}^q$ and $\vec{e}_j = (0,\cdots, 0, 1, 0, \cdots,0)
\in {\bf R}^q$ with $1$ being at the $j$-th entry.

Our first theorem shows that the Fisher information matrix of a GHMM can be written as an expectation similar to \eqref{Fisher-iid}.

\begin{theorem}\label{thm:Fisher} Assume conditions C1--C6 hold with $r = 2$. Then, we have
\begin{equation}\label{eqn:Fisher-rep}
I(\theta_0) = -E_{\omega_{\theta_0,2}}^{\theta_0}\left[ G(W_{1,\theta_0}^{(2)}, W_{0,\theta_0}^{(2)})\right],
\end{equation}
where $\omega_{\theta, r}$ is the stationary distribution of $\{((X_n, Y_n), W_{n,\theta}^{(r)}), n \geq 0\}$, $E_{\omega}^{\theta}$ is the
expectation taken when the above induced Markov chain is governed by $\theta$
and has an initial distribution equal to $\omega$, and
\begin{equation}\label{G}
    G(w_1, w_0) = \begin{pmatrix}
    g^{\nu(1,1)}(w_1, w_0) & \cdots & g^{\nu(1,q)}(w_1, w_0) \\
    \vdots & \ddots & \vdots \\
    g^{\nu(q,1)}(w_1, w_0) & \cdots & g^{\nu(q,q)}(w_1, w_0)
    \end{pmatrix},
\end{equation}
with $g^\nu$ defined in Lemma \ref{thm:additiv}, and for all $1 \leq j, k \leq q$,
\begin{equation*}
    \nu(j,k) = \vec{0} + \vec{e}_j + \vec{e}_k.
\end{equation*}
\end{theorem}

\begin{remark}
Note that one can link the function $G$ to the second derivatives of $\ell(\theta; Y_{0:1})$. See Remark \ref{remark:G-to-D2L} below for details.
\end{remark}

Our second theorem shows that the KL-divergence for GHMM can also be written in a form similar to \eqref{KL-iid}, and can be locally approximated by a quadratic
function determined by the Fisher information matrix as in \eqref{KL-converge}.

\begin{theorem}\label{thm:KL} Assume conditions C1--C2 hold with $r = 0$. Then, for any
$\theta_1 \in N_\delta(\theta_0)$, $K(\theta_1, \theta_0)$ is well-defined with
\begin{equation}\label{eqn:KLrep}
    K(\theta_1, \theta_0) = E_{\omega_{\theta_1,0}}^{\theta_1}\left[ g^0(W_{1,\theta_1}^{(0)}, W_{0,\theta_1}^{(0)})\right]
    -
    E_{\omega_{\theta_0,0}}^{\theta_1}\left[ g^0(W_{1,\theta_0}^{(0)}, W_{0,\theta_0}^{(0)})\right],
\end{equation}
where $E_{\omega_{\theta,0}}^{\theta_1}$ is defined as in Theorem \ref{thm:Fisher}, and function $g^0$ is defined in Lemma \ref{thm:additiv}. In addition, if the conditions C1--C6 hold with $r = 3$,
then as $\theta_1 \rightarrow \theta_0$, we have
\begin{equation}\label{HMMKL-converge}
K(\theta_1, \theta_0) = (\theta_1 - \theta_0)^t \frac{I(\theta_0)}{2} (\theta_1-\theta_0) + o(\Vert \theta_1 - \theta_0 \Vert^2).
\end{equation}
\end{theorem}

\begin{remark}
Note that one can link the function $g^0$ to $\ell(\theta; Y_{0:1})$. See Remark \ref{remark:KL-to-L} below for details.
\end{remark}

With the help of \eqref{eqn:Fisher-rep} and \eqref{eqn:KLrep}, we will prove the
following results related to the Cram\'{e}r-Rao lower bound and the H\'{a}jek-Le
Cam local asymptotic minimax theorem. The H\'{a}jek-Le Cam convolution theorem
for a finite state HMM can be found in  \cite{bickel1996inference}.
For any $n$, let $\mathcal{E}_n$ be the space of all estimators of $\theta$
based on $Y_{0:n}$, and $\mathcal{E}_n^U$ be the space of
all unbiased estimators  of $\theta$ based on $Y_{0:n}$. Denote $\hat{\theta}_n := \hat{\theta}_n(Y_{0:n})$
as an estimator based on $Y_{0:n}.$

\begin{theorem}\label{cor:CRLB} Assume conditions C1--C6 hold with $r = 2$.
Then, for any $v \in {\bf R}^q$ and $x \in \mathcal{X}$,
\begin{equation}\label{eqn:CRLB-classical}
	\lim_{n \rightarrow \infty} \inf_{\hat{\theta}_n \in \mathcal{E}_n^U} n E_x^{\theta_0}\left[ \left( v^t(\hat{\theta}_n(Y_{0:n}) - \theta_0) \right)^2 \right]   \geq v^t I^{-1}(\theta_0) v.
\end{equation}
In addition, assume C1--C6 hold with $r = 3$. Then, for $\delta = (nv^t I(\theta_0) v)^{-1/2}$, we have
\begin{align}\label{eqn:CRLB-general}
     \lim_{n \rightarrow \infty} \inf_{\hat{\theta}_n \in \mathcal{E}_n}  \max_{\theta \in \{\theta_0, \theta_0+\delta v \} } n E_x^{\theta}\left[ \Vert \hat{\theta}_n(Y_{0:n}) - \theta  \Vert^2 \right]
    \geq \frac{1}{16} \frac{\Vert v \Vert^2}{v^t I(\theta_0) v}.
\end{align}
\end{theorem}

\begin{remark}\label{remark:CRLB} In the case when $q=1$ (namely, $\theta \in
{\bf R}$), \eqref{eqn:CRLB-classical} reduces to
	\begin{equation}\label{eqn:CRLB-simple}
	\lim_{n \rightarrow \infty} \inf_{\hat{\theta}_n \in \mathcal{E}_n^U} n E_x^{\theta_0}\left[ \left( \hat{\theta}_n(Y_{0:n}) - \theta_0  \right)^2 \right]  \geq \frac{1}{I(\theta_0)}.
	\end{equation}
	In addition, for this one-dimensional case, one can generalize
	\eqref{eqn:CRLB-general} to
	\begin{align}\label{eqn:CRLB-1dim}
	\notag
	\lim_{c \rightarrow \infty} \lim_{n \rightarrow \infty} \inf_{\hat{\theta}_n \in \mathcal{E}_n} \sup_{\theta : \vert \theta - \theta_0 \vert
	\leq c/\sqrt{n}} n E_x^{\theta}\left[ \left( \hat{\theta}_n(Y_{0:n}) - \theta  \right)^2 \right]
	\geq  \frac{1}{16} \frac{1}{I(\theta_0)}.
	\end{align}
	See Chapter~8.7 of \cite{van2000asymptotic} for
	the execution on supreme  over a compact set.
\end{remark}

\begin{remark}\label{remark:CRLB-iid}
Equations~\eqref{eqn:CRLB-classical} and \eqref{eqn:CRLB-general} are similar to the
classical case where $Y_{0:n}$ are i.i.d. random variables.
In particular,  \eqref{eqn:CRLB-general} states
that as long as the estimator can shrink to a $n^{-1/2}$-neighborhood of
$\theta_0$, regardless of the constant term, then the square loss is uniformly
bounded from below. An interesting phenomenon here is that we still have the
same constant $\frac{1}{16}$ in \eqref{eqn:CRLB-general} as that in the i.i.d.\ case.
For this GHMM version, however, since we have no
characterization of the Fisher information matrix for fixed~$n$, the argument
requires that $n$ goes to infinity to link the mean square error to the Fisher
information $I(\theta_0)$.
\end{remark}

Finally, with the help of \eqref{eqn:KLrep}, we can prove the following
properties for KL-divergence in GHMM.
\begin{corollary}\label{cor:KL} Assume conditions C1--C2 hold with $r = 0$. Then, for any
$\theta_1 \in N_\delta(\theta_0)$,
\begin{enumerate}
    \item (Non-Negativity) $K(\theta_1, \theta_0) \geq 0$.
	 \item (Additivity) Suppose $Y_{0:n} = (Y_{0:n}^1, Y_{0:n}^2)$, and for any
	 $n \in \mathbb{N}$ and $y_{0:n} = (y_{0:n}^1, y_{0:n}^2)$, we have
	 $L(\theta; y_{0:n}) = L(\theta; y_{0:n}^1) L(\theta; y_{0:n}^2)$ for
	 $\theta = \theta_0, \theta_1$. Then,
    \begin{equation*}
        K(\theta_1, \theta_0) = K_1(\theta_1, \theta_0) + K_2(\theta_1, \theta_0),
    \end{equation*}
	 where for $i=1, 2$, $K_i$ is the KL-divergence defined by replacing
	 $L(\theta;Y_{0:n})$ in \eqref{KL-HMM} by $L(\theta;Y_{0:n}^i)$,
	 respectively.
\end{enumerate}
\end{corollary}

\begin{figure}[h]
    \centering
    \includegraphics[width=\linewidth]{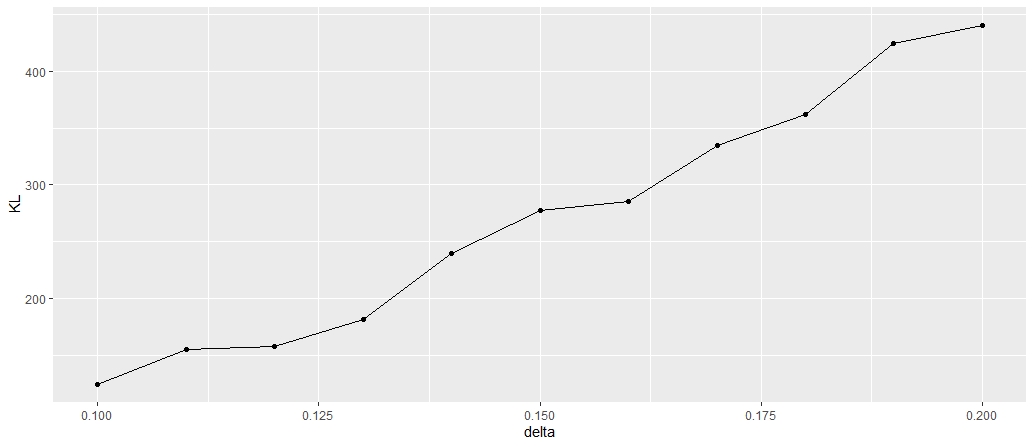}
      \caption{$K(\theta_1, \theta_0)$ for a three-state HMM}
      \label{fig:KL}
    \begin{flushleft}
	 {\footnotesize This figure presents $K(\theta_1, \theta_0)$ under different $\theta_1$. The model is defined in Remark~\ref{remark:convex}, where
	 $\theta_0$ corresponds to the model with $\delta=0$, and $\theta_1$ corresponds to the model with $\delta$ ranging from $0.1$ to $0.2$. For
	 each $\delta$, the $y$-axis represents the corresponding $K(\theta_1, \theta_0)$, which is computed via a Monte Carlo simulation using
	 \eqref{KL-HMM} with $X_0=0$ and $n=50,000$.}
    \end{flushleft}
\end{figure}

\begin{remark}\label{remark:convex}
Note that when $\{Y_n, n \geq 0\}$ is a sequence of i.i.d.\ finite mixture random variables or a Markov chain,
one can additionally prove that $K(\theta_1, \theta_0)$ is a convex function.
However, this is not the case here in general. To see why, let us consider the
case in which $\{Y_n, n \geq 0\}$ is a Markov chain. By using an argument similar
to Theorem~1 of \cite{rached2004kullback},  one can show that
\begin{align}\label{KL-MC}
\notag
    K(\theta_1, \theta_0) & = E_{\mu}^{\theta_1}\left[ \log f(Y_1; \theta_1|Y_0) \right] - E_{\mu}^{\theta_1}\left[ \log f(Y_1; \theta_0|Y_0) \right] \\
\notag
    & = E_{\mu}^{\theta_1}\left[ \log \frac{f(Y_1; \theta_1|Y_0) }{f(Y_1; \theta_0|Y_0) } \right] \\
    & = E_{\mu}^{\theta_1}\left[ \log \frac{f(Y_1; \theta_1|Y_0)\mu(Y_0) }{f(Y_1; \theta_0|Y_0)\mu(Y_0)} \right],
\end{align}
where $\mu(\cdot)$ is the invariant measure of $\{Y_n, n \geq 0\}$ under $\theta_1$, and
$E_{\mu}^{\theta_1}$ is the expectation when the Markov chain is
governed by $P^{\theta_1}$ and $Y_0$ is $\mu$-distributed. By such,
the classical argument applying log-sum inequality
leads to the convexity of $K(\theta_1,
\theta_0)$.

This argument, however, does not work for the case where $\{Y_n, n \geq 0\}$ is
a HMM. This is because unlike \eqref{KL-MC}, the two expectations in
\eqref{eqn:KLrep} are under different invariant measures, so they cannot be
combined as in \eqref{KL-MC}, and the log-sum inequality cannot be
applied. This non-convexity becomes a unique feature for HMM that is different
from an i.i.d.\ or Markov chain scenario. Similar non-convexity for the KL-divergence is observed in \cite{vidyasagar2010kullback}.

We end this remark by a numerical illustration. Consider a three-state HMM with
$\mathcal{X} = \{1, 2, 3\}$, for which $P\{X_1 = x_1 | X_0 = x_0 \} = \frac{1}{3}$
for all $x_0, x_1 \in \mathcal{X}$. As for the observations, we assume $Y_n \in \{1,
2\}$ with $P\{ Y_n = 1 | X_n = x\} = q_x^\delta$, where
\begin{equation*}
(q_1^\delta, q_2^\delta, q_3^\delta) = \left( 1, \frac{1}{2}+\delta, 0 \right).
\end{equation*}
Let $\theta_0$ be  the corresponding probability measure with $\delta=0$, and
$\theta_1$ be the corresponding  probability measure with $\delta$ ranging from
$0.1$ to $0.2$. The computed $K(\theta_1, \theta_0)$ is presented in
Figure~\ref{fig:KL}. As expected, $K(\theta_1, \theta_0)$ is decreasing when $\delta$ decreases. However, the
figure shows that $K(\theta_1, \theta_0)$ is not convex, as mentioned above.
\end{remark}

\begin{remark}
The conditions can be slightly relaxed. For example, as
Theorems~\ref{thm:Fisher}--\ref{cor:CRLB} and Corollary~\ref{cor:KL} only involve the
neighborhood of $\theta_0$, the differentiability assumption in C2 can be
relaxed to only $\theta \in N_\delta(\theta_0)$ instead of all $\theta$. Also,
as one can see in the proof, the second part of C4 is used only to bound the
third-order derivatives of $\log L(\theta; Y_{0:n})$ in order to obtain the
small-$o$ term in \eqref{HMMKL-converge}. By such, C4 can be relaxed for the
results unrelated to the small-$o$ term in \eqref{HMMKL-converge}.
\end{remark}

\section{AIC Model Selection}\label{sec:AIC}

In this section, we will use Akaike's information criterion (AIC) in determining the order of a GHMM.
Note that the HMM is defined in a general sense as that in Section \ref{sec:intro}.
To this end, we present an objective procedure for the determination of the order of an ergodic general hidden Markov model with a finite state space. The procedure exploits the asymptotic properties
of the likelihood ratios statistics in \cite{fuh_pang_2022}, and the KL-divergence, defined in Theorem \ref{thm:KL}, for the discrimination between two GHMM distributions.

 $\{X_n, n\ge 0\}$ is called a $k$-order Markov chain if $k$ is the smallest non-negative integer such that
\begin{eqnarray*}
	P(X_n|X_{n-1}, X_{n-2},\cdots)=P(X_n|X_{n-1},X_{n-2},\cdots,X_{n-k})\quad {\rm for~all~}n.
\end{eqnarray*}
In what follows we assume that $\{X_n, n\ge 0\}$ is an $m$-order Markov chain on a finite state space $\mathcal{D}=\{1,\cdots,l \}$. It is known that for an $m$-order Markov chain  $\{X_n, n\ge 0\}$, then $\{(X_t,X_{t+1},\cdots,X_{t+m-1}), t=0,\dots,n-m+1,\cdots \}$ forms
a Markov chain.
Following the definition in Section \ref{sec:intro} and
(\ref{eqn:lik}), $\{Y_n, n\ge 0\}$ is called an $m$-order GHMM. In what follows in this section, we further assume conditions C1-C6 hold for this $m$-order GHMM.

Let $Y_{0:n}=\{Y_0,\cdots,Y_n\}$ be the observations from an $m$-order GHMM such that $Y_t \in {\bf R}^d$ with a distribution depending on $X_n$ and $Y_{n-1}$, but independent to others. Let $f(\cdot; {_m\theta}|x, y)$ be the pdf of $Y_n$ given $X_n=x$ and $Y_{n-1}=y$, with respect to a $\sigma$-finite measure on ${\bf R}^d$, where ${_m\theta}=({_m\theta_1},\cdots,{_m\theta_q})^t\in \Theta\subset {\bf R}^q$ with $q\ge m$. Furthermore, let $f(\cdot; {_m\theta}|x_0)$ be the pdf of $Y_0$ given $X_0 = x_0$.
Denote $i_j=(x_j,\cdots,x_{j+m-1})$ for $j=0,\cdots, n-m+1$, then the full likelihood of this $m$-order GHMM is
\begin{align}
	L({_m\theta}; Y_{0:n})
	= \sum_{x_0=1}^l \cdots  \sum_{x_n=1}^l \pi_{_m\theta}(i_0) f(Y_0;{_m\theta}|x_0)\prod_{j=1}^{n-m+1} p_{_m\theta}(i_{j-1},i_j) \prod_{j=1}^{n}
 f(Y_j;{_m\theta}| x_j, Y_{j-1}),
\end{align}
where $\pi_{_m\theta}(\cdot)$ is the stationary distribution of the  Markov chain  $\{(X_t,X_{t+1},\cdots,X_{t+m-1}), t=0,\dots,n-m+1,\cdots \}$  and $p_{_m\theta}(\cdot,\cdot)$ is the corresponding transition probability kernel. Denote
\begin{align*}
	_m\hat{\theta}= \arg\max_{_m\theta\in \Theta}L({_m\theta}; Y_{0:n}).
\end{align*}
That is, ${_m\hat{\theta}}$ is the maximum likelihood estimator (MLE) of ${_m\theta}$ based on the observations $Y_{0:n}$ from this $m$-order GHMM. In what follows, we suppose the true value of ${_m\theta}$ is ${_m\ddot{\theta}}$. Denote $I({_m\ddot{\theta}})$ as the Fisher information matrix which corresponds to this $m$-order GHMM. Then, for any $\theta \in \Theta\subset {\bf R}^q$ we denote
\[\|\theta\|_J^2=\theta^tI({_m\ddot{\theta}})\theta.\]

Let $\Delta$ be the difference operator on the superscript, $\Delta_t^j=t^j-t^{j-1}$ for
$ j\ge 1.$
As in \cite{akaike1973information}, we consider the case where $k\le m$. Suppose that ${_k{\theta}}$ is restricted to the parameter space $\Theta$ with $\Delta_l^{m+1}-\Delta_l^{k+1}$ components of ${_k{\theta}}$ being equal to zero due to the change of order from $m$ to $k$ for the Markov chain $\{X_n, n\ge 0\}$.
Without loss of generality, we suppose that the last $\Delta_l^{m+1}-\Delta_l^{k+1}$ components of ${_k\theta}$ are restricted to be zero when the order of the Markov chain $\{X_n, n\ge 0\}$ changes from $m$ to $k$. Denote the above restricted parameter space as $\Theta_R$, and let ${_k\hat{\theta}}$ be the MLE of ${_k{\theta}}$ in this restricted parameter space. Moreover, let
${_k\ddot{\theta}}\in \Theta_R$ such that
\begin{align*}
	\|{_k\ddot{\theta}}-{_m\ddot{\theta}}\|^2_J=\min_{\theta\in \Theta_R}\|\theta-{_m\ddot{\theta}}\|^2_J.
\end{align*}
That is, ${_k\ddot{\theta}}$ is the projection of ${_m\ddot{\theta}}$ in the space $\Theta_R$ with respect to the metrics defined by $\|\cdot\|_J$. Denote $p=q-(\Delta_l^{m+1}-\Delta_l^{k+1})$, which is the active dimension of the restricted parameter space $\Theta_R$. Then, it is easy to see that
\begin{align*}
	\sum_{j=1}^p I({_m\ddot{\theta}})_{i,j}\times {_k\ddot{\theta}_{j}}=\sum_{j=1}^q I({_m\ddot{\theta}})_{i,j}\times {_k\ddot{\theta}_{j}}=\sum_{j=1}^q I({_m\ddot{\theta}})_{i,j}\times {_m\ddot{\theta}_{j}},\quad {\rm for~all}~i=1,\cdots,q,
\end{align*}
where $I({_m\ddot{\theta}})_{i,j}$ stands for the $(i,j)$-th element of $I({_m\ddot{\theta}})$. This further implies that
\begin{align}\label{projection}
	\|{_k\hat{\theta}}-{_m\ddot{\theta}}\|^2_J=\|{_k\hat{\theta}}-{_k\ddot{\theta}}\|^2_J+\|{_k\ddot{\theta}}-{_m\ddot{\theta}}\|^2_J.
\end{align}
In what follows, we take $\|{_k\hat{\theta}}-{_m\ddot{\theta}}\|^2_J$ as the loss function, because it is approximately equal to $2K({_m\ddot{\theta}}, {_k\hat{{\theta}}})$, which is close to $2K({_m\ddot{\theta}}, {_k{\theta}})$ with $K({_m\ddot{\theta}}, {_k{\theta}})$ denoting the KL-divergence between the $m$-order GHMM and the $k$-order GHMM.

Denote
\[{_k\lambda_m}=\frac{L({_m\hat{\theta}}; Y_{0:n})}{L({_k\hat{\theta}}; Y_{0:n})}\]
as the ratio of the maximum likelihood given that $Y_{0:n}$ is from an $m$-order GHMM to that given that $Y_{0:n}$ is from a $k$-order GHMM. Then
\begin{align}\label{ratio}
	\log ({_k\lambda_m})=&\log L({_m\hat{\theta}}; Y_{0:n})-\log L({_k\hat{\theta}}; Y_{0:n}) \notag\\
	=& \log \frac{L({_k\ddot{\theta}}; Y_{0:n})}{L({_k\hat{\theta}}; Y_{0:n})}-\log \frac{L({_k\ddot{\theta}}; Y_{0:n})}{L({_m\hat{\theta}}; Y_{0:n})}.
\end{align}
Taking into account the relations
\begin{align*}
	\frac{\partial \log L({_m\theta}; Y_{0:n})}{\partial {_m\theta}}{\Big|}_{{_m\theta}={_m\hat{\theta}}}=0 \quad \mbox{and}\quad \frac{\partial \log L({_k\theta}; Y_{0:n})}{\partial {_k\theta}}{\Big|}_{{_k\theta}={_k\hat{\theta}}}=0,
\end{align*}
it follows from the Taylor expansions that
\begin{align}\label{Taylor1}
	&\log L({_k\ddot{\theta}}; Y_{0:n})=\log L({_m\hat{\theta}}; Y_{0:n}) \notag\\
	&\quad +\frac{n}{2}\sum_{i=1}^q\sum_{j=1}^q ({_k\ddot{\theta}_i}-{_m\hat{\theta}_i})({_k\ddot{\theta}_j}-{_m\hat{\theta}_j}) \frac{1}{n}\frac{\partial^2 \log L({_m\theta}; Y_{0:n})}{\partial {_m\theta_i}\partial {_m\theta_j}}{\Big|}_{{_m\theta}={_m\hat{\theta}}+\alpha({_k\ddot{\theta}}-{_m\hat{\theta}})}
\end{align}
with $0\le \alpha\le 1$, and
\begin{align*}
	&\log L({_k\ddot{\theta}}; Y_{0:n})=\log L({_k\hat{\theta}}; Y_{0:n}) \notag\\
	&\quad +\frac{n}{2}\sum_{i=1}^p\sum_{j=1}^p ({_k\ddot{\theta}_i}-{_k\hat{\theta}_i})({_k\ddot{\theta}_j}-{_k\hat{\theta}_j})\frac{1}{n}\frac{\partial^2 \log L({_k\theta}; Y_{0:n})}{\partial {_k\theta_i}\partial {_k\theta_j}}{\Big|}_{{_k\theta}={_k\hat{\theta}}+\beta({_k\ddot{\theta}}-{_k\hat{\theta}})}\nonumber
\end{align*}
with $0\le \beta\le 1$. Applying the decomposition in (\ref{function_g}), we write
\begin{align*}
	\log L({_m\theta}; Y_{0:n})=\sum_{t=1}^n g^0( M_t^0({_m\theta}), M_{t-1}^0({_m\theta})) +g_0^0(M_0^0({_m\theta})),
\end{align*}
where the definitions of $g^0$ and $g_0^0$ can be found in (\ref{function_gg}). Then, we have
\begin{align*}
	&\frac{1}{n}\frac{\partial^2 \log L({_m\theta}; Y_{0:n})}{\partial {_m\theta_i}\partial {_m\theta_j}}{\Big|}_{{_m\theta}={_m\hat{\theta}}+\alpha({_k\ddot{\theta}}-{_m\hat{\theta}})}\\
	=& \frac{1}{n}\sum_{t=1}^n \frac{\partial^2 g^0( M_t^0({_m\theta}), M_{t-1}^0({_m\theta}))}{\partial {_m\theta_i}\partial {_m\theta_j}}{\Big|}_{{_m\theta}={_m\hat{\theta}}+\alpha({_k\ddot{\theta}}-{_m\hat{\theta}})} +\frac{1}{n}\frac{\partial^2 g_0^0(M_0^0({_m\theta}))}{\partial {_m\theta_i}\partial {_m\theta_j}}{\Big|}_{{_m\theta}={_m\hat{\theta}}+\alpha({_k\ddot{\theta}}-{_m\hat{\theta}})}\\
	\rightarrow & -I({_m\ddot{\theta}})_{ij} ~~~P^{_m\ddot{\theta}} \mbox{-a.s.}
\end{align*}
provided $\sqrt{n} \|_k\ddot{\theta} - {_m\theta}\|_J$ is bounded since ${_m\hat{\theta}}$ actually is the MLE of ${_m\ddot{\theta}}$ which is asymptotically efficient (cf. Theoreom 2 in \cite{fuh_pang_2022}). This together with (\ref{Taylor1}) imply that
\begin{align*}
	\log \frac{L({_k\ddot{\theta}}; Y_{0:n})}{L({_m\hat{\theta}}; Y_{0:n})}=\frac{n}{2}({_k\ddot{\theta}}-{_m\hat{\theta}})I({_m\ddot{\theta}})({_k\ddot{\theta}}-{_m\hat{\theta}})^t+o_{P^{_m\ddot{\theta}}}(1)=\frac{n}{2}\|{_k\ddot{\theta}}-{_m\hat{\theta}}\|^2_J+o_{P^{_m\ddot{\theta}}}(1).
\end{align*}
Similarly, it can be proved that if $\sqrt{n}\|{_k\ddot{\theta}}-{_m\theta}\|_J$ is bounded,
\begin{align*}
	\log \frac{L({_k\ddot{\theta}}; Y_{0:n})}{L({_k\hat{\theta}}; Y_{0:n})}=\frac{n}{2}({_k\ddot{\theta}}-{_k\hat{\theta}})I({_m\ddot{\theta}})({_k\ddot{\theta}}-{_k\hat{\theta}})^t+o_{P^{_m\ddot{\theta}}}(1)=\frac{n}{2}\|{_k\ddot{\theta}}-{_k\hat{\theta}}\|^2_J+o_{P^{_m\ddot{\theta}}}(1).
\end{align*}
Thus, it follows from (\ref{ratio}) that
\begin{align}\label{decomposition}
	& -2\log ({_k\lambda_L}) \nonumber\\
=& n\|{_k\ddot{\theta}}-{_m\hat{\theta}}\|^2_J-n\|{_k\ddot{\theta}}-{_k\hat{\theta}}\|^2_J+o_{P^{_m\ddot{\theta}}}(1)\nonumber\\
=& n\|{_k\ddot{\theta}}-{_m\ddot{\theta}}\|^2_J+ n\|{_m\ddot{\theta}}-{_m\hat{\theta}}\|^2_J-n\|{_k\ddot{\theta}}-{_k\hat{\theta}}\|^2_J -2n({_k\ddot{\theta}}-{_m\ddot{\theta}},{_m\hat{\theta}}-{_m\ddot{\theta}})_J+o_{P^{_m\ddot{\theta}}}(1),
\end{align}
where $({_k\ddot{\theta}}-{_m\ddot{\theta}},{_m\hat{\theta}}-{_m\ddot{\theta}})_J=({_k\ddot{\theta}}-{_m\ddot{\theta}})^t I({_m\ddot{\theta}})({_m\hat{\theta}}-{_m\ddot{\theta}})$ denotes the inner product of $({_k\ddot{\theta}}-{_m\ddot{\theta}})$ and $({_m\hat{\theta}}-{_m\ddot{\theta}})$ defined by the Fisher information matrix $I({_m\ddot{\theta}})$.

By Theorem 2 in \cite{fuh_pang_2022}, we have
\begin{align*}
	n\|{_m\hat{\theta}}-{_m\ddot{\theta}}\|^2_J\rightarrow \chi^2_q(0)~~{\rm in~distribution}.
\end{align*}
Note that geometrically ${_k\hat{\theta}} - {_k\ddot{\theta}}$ is approximately the projection of ${_m\hat{\theta}}-{_m\ddot{\theta}}$ into the space of $\Theta_R$, therefore as long as $\sqrt{n}\|{_k\ddot{\theta}}-{_m\ddot{\theta}}\|_J$ is bounded, it is true that
\begin{align*}
	n\|{_m\hat{\theta}}-{_m\ddot{\theta}}\|^2_J-n\|{_k\hat{\theta}}-{_k\ddot{\theta}}\|^2_J &\rightarrow \chi^2_{q-p}(0)~~{\rm in~distribution},\\
	n\|{_k\hat{\theta}}-{_k\ddot{\theta}}\|^2_J &\rightarrow \chi^2_{p}(0)~~{\rm in~distribution},
\end{align*}
and $n\|{_m\hat{\theta}}-{_m\ddot{\theta}}\|^2_J-n\|{_k\hat{\theta}}-{_k\ddot{\theta}}\|^2_J$ and $n\|{_k\hat{\theta}}-{_k\ddot{\theta}}\|^2_J$ are asymptotically independent. Note that Theorem 2 of \cite{fuh_pang_2022} also implies that the standard deviation of the asymptotic distribution of $n({_k\ddot{\theta}}-{_m\ddot{\theta}},{_m\hat{\theta}}-{_m\ddot{\theta}})_J$ is equal to $\sqrt{n}\|{_k\ddot{\theta}}-{_m\ddot{\theta}}\|_J$. Thus, $n({_k\ddot{\theta}}-{_m\ddot{\theta}},{_m\hat{\theta}}-{_m\ddot{\theta}})_J$ is negligible in comparison with the term $n\|{_k\ddot{\theta}}-{_m\ddot{\theta}}\|^2_J$ if the latter is large enough.

Then, it follows from the above arguments and the equations (\ref{projection}) and (\ref{decomposition}) that $\{p+[-2\log ({_k\lambda_m})-(q-p)]\}/n=[-2\log ({_k\lambda_m})-q+2p]/n$ serves as a useful estimator of $E^{{_m\ddot{\theta}}}\|{_k\hat{\theta}}-{_m\ddot{\theta}}\|^2_J$. That is, the determining the order of a GHMM can be conducted via minimizing the following AIC criterion
\begin{align}\label{eqn:criterion}
-2\log ({_k\lambda_m})-q+2p = -2\log ({_k\lambda_m})-2(\Delta_l^{m+1}-\Delta_l^{k+1})+q
\end{align}
by recalling that $p=q-(\Delta_l^{m+1}-\Delta_l^{k+1})$. Getting rid of the terms which are independent of $k$ from the RHS of (\ref{eqn:criterion}), we arrive at the following AIC criterion for GHMM:
\begin{align}\label{eqn:criterion1}
{\rm AIC}(k)=2\log L({_k\hat{\theta}}; Y_{0:n})+2\Delta_l^{k+1}.
\end{align}

By making use of the same method, it can be shown that for a given $m$-order GHMM, the AIC criterion for selecting the number of hidden states is
\begin{align}\label{eqn:criterion2}
{\rm AIC}(k)=2\log L({_k\hat{\theta}}; Y_{0:n})+2\Delta_k^{m+1},
\end{align}
where $k$ denotes the number of hidden states in the $m$-order GHMM, and ${_k\hat{\theta}}$ denotes the MLE of the parameter in such model. \cite{Yonekuraa_Beskosa_Singhb:2021} considers AIC model selection for HMM when $Y_n$ depends on $X_n$ only.

\section{Examples}\label{sec:Examples}

\begin{example}\label{Ex1} Gaussian Vector Autoregressive-Moving Average Model.

Consider a Gaussian vector autoregressive-moving average (VARMA) model of order $(p,q)$ in $m$ dimension  (see \cite{klein2008asymptotic}) such that, for all $n \geq 0$,
\begin{equation}\label{VARMA}
\sum_{j=0}^p \alpha_j Y_{n-j} = \sum_{j=0}^q \beta_j \epsilon_{n-j},
\end{equation}
in which $\alpha_j$ and $\beta_j$ are $m$-by-$m$ real-valued matrices with $\alpha_0 = \beta_0 = I_m$ (the $m$-by-$m$ identity matrix), and $\{\epsilon_n, n \geq 1\}$ are i.i.d. $m$-dimensional normal random variables with zero mean and covariance matrix $\Sigma$. Here, we assume that $p$, $q$ and $\Sigma$ are known, so the unknown parameter can be denoted as
\begin{equation*}
\theta = vec\left\{ \alpha_1, \cdots, \alpha_p, \beta_1, \cdots, \beta_q \right\},
\end{equation*}
where for any matrix $M$, $vec \left\{ M \right\}$ denotes the vector created by stacking all the columns of $M$ on top of each other.

It is known that the VARMA model in \eqref{VARMA} can be represented as a linear state space model (LSSM) in various ways, cf. \cite{harvey1979maximum}. For example, let $h = \max(p,q)$, and suppose $X_n \in {\bf R}^{hm}$ and $Y_n \in {\bf R}^m$ satisfying the LSSM
\begin{align}\label{LSSM}
\notag
X_{n+1} & = \Phi X_n + F \epsilon_n, \\
Y_n & = H X_n + \epsilon_n,
\end{align}
where
\begin{equation*}
\Phi =
\begin{pmatrix}
-\alpha_1 & I_m & 0_m & \cdots & 0_m \\
-\alpha_2 & 0_m & I_m & \cdots & 0_m \\
\vdots & \vdots & \vdots & \ddots & \vdots \\
-\alpha_{h-1} & 0_m & 0_m & \cdots & I_m \\
-\alpha_h & 0_m & 0_m & \cdots & 0_m
\end{pmatrix}_{hm \times hm}, ~~~
F =
\begin{pmatrix}
\beta_1 - \alpha_1 \\
\beta_2 - \alpha_2 \\
\vdots \\
\beta_h - \alpha_h
\end{pmatrix}_{hm \times m},
\end{equation*}
and $H = (I_m, 0_m, \cdots, 0_m)$, with $0_m$ being the $m$-by-$m$ zero matrix, $\alpha_i = 0_m$ for all $i >p$ and $\beta_j = 0_m$ for all $j > q$. Then, \cite{harvey1979maximum} has shown that the $Y_n$ in \eqref{LSSM} satisfies the VARMA model in \eqref{VARMA}; see also \cite{klein2008asymptotic}, \cite{pearlman1980algorithm} and \cite{hannan2012statistical}.

Since \eqref{LSSM} is in the form of LSSM, consider the sample innovation $\hat{\epsilon}_n$ obtained by the following Kalman filter equations:

\begin{align}\label{Kalman}
\begin{cases}
P_{n+1}  = \Phi P_n \Phi^t + \Sigma - (\Phi P_n H^t)(H P_n H^t)^{-1}(H P_n \Phi^t), \\
K_n  = (\Phi P_n H^t) (H P_n H^t)^{-1}, \\
\hat{X}_{n+1|n}  = (\Phi - K_n H) \hat{X}_{n|n-1} + K_n Y_n, \\
\hat{Y}_{n|n-1}  = H \hat{X}_{n|n-1}, \\
\hat{\epsilon}_n  = Y_n - \hat{Y}_{n|n-1},
\end{cases}
\end{align}
with $P_1$ given by $P_1=\Phi P_1 \Phi^t + F F^t$. Further denote $\hat{\Sigma}_n = E\left[ \hat{\epsilon}_n \hat{\epsilon}_n^t \right]$. Then, the log likelihood of the VARMA can be written as
		\begin{equation*}\label{logLik-VARMA}
		\ell(\theta;Y_{1:n}) = \sum_{i=1}^n \left\{
		-\frac{n}{2} \log 2\pi - \frac{1}{2} \log |\hat{\Sigma}_i| - \frac{1}{2} \hat{\epsilon}_i^t \hat{\Sigma}_i^{-1} \hat{\epsilon}_i \right\},
		\end{equation*}
		which further implies that
		\begin{align}\label{2nd-VARMA}
		\notag
		& - \frac{1}{n} E^{\theta_0}\left[ \frac{\partial^2 \ell(\theta_0; Y_{1:n})}{\partial \theta \partial \theta^t} \right] \\
		\notag
		= & \frac{1}{n} \sum_{i=1}^n \frac{1}{2} \left( \frac{\partial vec \left\{ \hat{\Sigma}_i \right\} }{\partial \theta}\right)^t (\hat{\Sigma}_i \otimes \hat{\Sigma}_i)^{-1} \left( \frac{\partial vec \left\{ \hat{\Sigma}_i \right\} }{\partial \theta}\right) \\
		& +
		\frac{1}{n} \sum_{i=1}^n E^{\theta_0}\left[ \left(\frac{\partial \hat{\epsilon}_i}{\partial \theta} \right)^t \hat{\Sigma}_i^{-1} \left(\frac{\partial \hat{\epsilon}_i}{\partial \theta} \right) \right],
		\end{align}
		where $\otimes$ denotes the Kronecker product. See \cite{klein2008asymptotic} for details.

 To apply Theorems \ref{thm:Fisher} and \ref{thm:KL}, we need to check conditions C1-C6 hold. Section 6.2 of \cite{Fuh2006} checks that conditions C1-C6 hold for ARMA. As for LSSM, Section 16.5.1. of \cite{Meyn&Tweedie} shows that C1 (the $\omega$-uniformity) hold for any dimensional LSSM, therefore by using the same argument, C1 holds for VARMA (\ref{VARMA}) and (\ref{LSSM}). By using the normality of $\epsilon_n$, it is straightforward to check conditions C2-C6 hold.

Now, \cite{klein2008asymptotic} has shown that, when $n \rightarrow \infty$, the limiting distribution of $\hat{\epsilon}_n$ is the same as the distribution of $\epsilon_1$. In addition, since both $\hat{\epsilon}_n$ and $\epsilon_1$ are Gaussian, the covariance matrix of $\hat{\epsilon}_n$ converges to that of $\epsilon_1$; in other words, $\hat{\Sigma}_n \rightarrow \Sigma$, which is independent to $\theta$.
By such, the first term in \eqref{2nd-VARMA} goes to zero as $n \rightarrow \infty$, and therefore,
		\begin{align}\label{VARMA-3}
		I(\theta_0) = \lim_{n \rightarrow \infty} \frac{1}{n} \sum_{i=1}^n E^{\theta_0}\left[ \left(\frac{\partial \hat{\epsilon}_i}{\partial \theta} \right)^t  \Sigma^{-1} \left(\frac{\partial \hat{\epsilon}_i}{\partial \theta} \right) \right].
		\end{align}
In addition, it is known that $K_n \rightarrow K_\infty$ as $n \rightarrow \infty$ for some finite constant matrix $K_\infty$ (depending on $\theta$). Therefore the asymptotic version of \eqref{Kalman} becomes
\begin{align}\label{Kalman_asymp}
\begin{cases}
\hat{X}_{n+1|n}^\infty  = (\Phi - K_\infty H) \hat{X}_{n|n-1}^\infty   + K_\infty Y_n, \\
\hat{Y}_{n|n-1}^\infty    = H \hat{X}_{n|n-1}^\infty  , \\
\hat{\epsilon}_n^\infty    = Y_n - \hat{Y}_{n|n-1}^\infty.
\end{cases}
\end{align}
Then by \eqref{VARMA-3}, we have
\begin{align}\label{VARMA-3-infty}
I(\theta_0) = \lim_{n \rightarrow \infty} \frac{1}{n} \sum_{i=1}^n E^{\theta_0}\left[ \left(\frac{\partial \hat{\epsilon}_i^\infty}{\partial \theta} \right)^t  \Sigma^{-1} \left(\frac{\partial \hat{\epsilon}_i^\infty}{\partial \theta} \right) \right].
\end{align}

Now, by \eqref{Kalman_asymp}, we have
\begin{align}
    \frac{\partial \hat{X}_{n+1|n}^\infty}{\partial \theta}
    =
    \frac{\partial (\Phi - K_\infty H)}{\partial \theta}\hat{X}_{n+1|n}^\infty
    + (\Phi - K_\infty H) \frac{\partial \hat{X}_{n+1|n}^\infty}{\partial\theta} + \frac{\partial K_\infty}{\partial \theta} Y_n,
\end{align}
which indicates that $Z_n = \left( X_n, Y_n, \hat{X}_{n|n-1}, \frac{\partial}{\partial \theta} \hat{X}_{n|n-1} \Big\vert_{\theta = \theta_0} \right)$ forms a Markov chain. Moreover, note that
$ \hat{\epsilon}_n^\infty  = Y_n - \hat{Y}_{n|n-1}^\infty = Y_n - H\hat{X}_{n|n-1}^\infty$,
which~implies
\begin{align}
 \frac{\partial \hat{\epsilon}_n^\infty}{\partial \theta} & =  - \frac{\partial H}{\partial \theta}\hat{X}_{n|n-1}^\infty - H\frac{\partial \hat{X}_{n|n-1}^\infty}{\partial \theta}.
\end{align}
In other words, $\frac{\partial}{\partial \theta}\hat{\epsilon}_n^\infty$ is a function of $Z_n$. Therefore, Theorem \ref{thm:Fisher} essentially shows that $Z_n$ is stationary under $\omega_{\theta_0,2}$, and therefore
\begin{align}\label{LSSMF}
I(\theta_0)  =  E_{\omega_{\theta_0,2}}^{\theta_0}\left[ \left( \frac{\partial \hat{\epsilon}_1^\infty}{\partial \theta} \right)^t \Sigma^{-1} \left( \frac{\partial \hat{\epsilon}_1^\infty}{\partial \theta}\right){\Big|}_{\theta=\theta_0} \right].
\end{align}
The result (\ref{LSSMF}) is the same as that in \cite{klein2008asymptotic}, in which they use $\frac{\partial \epsilon}{\partial \theta}$ to denote a random variable with distribution as the limiting distribution of $\frac{\partial \hat{\epsilon}_n}{\partial \theta}$ as $n \rightarrow \infty$.
Here we derive the Fisher information (\ref{LSSMF}) from the
invariant probability measure of the enlarged Markov chain point of view.
The Fisher information for a general LSSM can be found in
\cite{klein2000direct}.

Note that the result above is made possible due to the fact that, if $(X_n, Y_n)$ follows the LSSM with parameter $\theta$, then under $P^{\theta}$, the limiting distribution of $\hat{\epsilon}_n$ is the same as $\epsilon_1$. If we want to find a representation of the KL-divergence, then we need to evaluate the limiting distribution of $\hat{\epsilon}_n$ under $P^{\theta'}$ for $\theta' \neq \theta$. This is due to that the second expectation in (\ref{eqn:KLrep}) involves $W_{1,\theta_0}^{(0)}$ and $W_{0,\theta_0}^{(0)}$ under $P^{\theta_1}$ with $\theta_0 \neq \theta_1$. Instead,
(\ref{HMMKL-converge}) in Theorem \ref{thm:KL}
provides a local approximation of the KL-divergence in terms of Fisher information.
Moreover, we provide a theoretical justification of using AIC model selection criterion in (\ref{eqn:criterion2}) to choose $(p,q)$ in (\ref{VARMA}) and (\ref{LSSM}). A computational method
of the KL-divergence and  AIC model selection for LSSM
can be found in \cite{Bengtsson_Cavanaugh_2006}.
\end{example}

\begin{example}\label{Ex2}
Recurrent Neural Networks.

To start with, we consider the following linear recurrent neural network (RNN) as
\begin{eqnarray}\label{6.8a}
Y_n = \mu_{y,n} + \sigma_{y,n}  \varepsilon_n,
\end{eqnarray}
where $(\mu_{y,n}, \sigma^2_{y,n})\sim \varphi_{\tau}(X_{n-1})$, with $\varphi_{\tau}$ can be any highly flexible function such as neural networks. $\sigma_{y,n} > 0$ P-a.s., $\varepsilon_n \sim N(0,1)$ is a sequence of i.i.d. random variables, and $\varepsilon_n$ is independent of $\{Y_{n-k}, k \geq 1\}$ for all $n$.

To illustrate the GHMM approach for the  linear RNN.
By using (\ref{6.8a}) as the output model for $Y_n$,
the linear RNN updates its hidden state using the recurrence equation:
\begin{eqnarray}\label{RNN5a}
X_n = f_\theta( \phi_{\tau}(Y_n),  X_{n-1}) = \delta +  \alpha Y_{n-1} +  \beta	X_{n-1},
\end{eqnarray}
where $\delta > 0,~ \alpha > 0,$  and $\beta > 0$ are constants.

As noted in \cite{fuh2021a} that the linear RNN model (\ref{6.8a}) and (\ref{RNN5a})  can be regarded as the celebrated GARCH$(1,1)$ model
when $\mu_{y,n} = 0$ and $\sigma_{y,n}^2 = X_{n}$ in (\ref{6.8a}), and $Y_{n-1}$ is replaced by $Y_{n-1}^2 $ in  (\ref{RNN5a}).
However, $\mu_{y,n}$ and $\sigma_{y,n}$, defined in \eqref{6.8a}, can be nonlinear functions of $X_{n}$ in
general.

To have an explicit computation of the Fisher information and KL-divergence, we
consider a specific form of (\ref{6.8a}),  a simple GARCH$(1,1)$ model (see \cite{bollerslev1986generalized}), as follows:
	\begin{equation*}\label{GARCH}
	Y_n = \sigma_n \epsilon_n, ~~~ \mbox{and} ~ \sigma_n^2 = \delta +  \alpha Y_{n-1}^2 + \beta \sigma_{n-1}^2,
	\end{equation*}
	where $\delta > 0$, $\alpha > 0$ and $\beta > 0$ are constants with $\alpha+\beta<1$, and $\{\epsilon_n, n \geq 1\}$ are i.i.d. standard normal random variables with $\epsilon_n$ independent of $\{Y_t, t=1, \cdots, n-1\}$.

Section 6.3 of \cite{Fuh2006} has checked that conditions C1-C6 hold for the GARCH$(p,q)$ model. As for the linear RNN, condition C1 (the $\omega$-uniformity) can be found in Section 16.5.1 of \cite{Meyn&Tweedie} using the state space representation. Conditions C2-C6 hold due to the normality assumption of $\varepsilon_n$. Therefore Theorems \ref{thm:Fisher} and \ref{thm:KL} can be applied.
	
	Denote $\theta=(\delta, \alpha, \beta)^t$. Note that the log likelihood function of GARCH$(1,1)$ model can be expressed as
	\begin{equation}\label{GARCH-sum}
	\ell(\theta; Y_{1:n}) = \sum_{t=1}^n \left\{ -\frac{1}{2}\log 2\pi - \frac{1}{2} \log \sigma_t^2 - \frac{1}{2} \frac{Y_t^2}{\sigma_t^2}\right\},
	\end{equation}
	which further implies
	\begin{equation}\label{GARCH-Fisher-0}
	-E^{\theta_0}\left[ \frac{\partial^2 \ell(\theta; Y_{1:n})}{\partial \theta \partial \theta^t} \right]
	= \sum_{t=1}^n E^{\theta_0}\left[ \frac{1}{2\sigma_t^2} \frac{\partial \sigma_t^2}{\partial \theta}\frac{\partial \sigma_t^2}{\partial \theta^t} \right].
	\end{equation}
	See \cite{ma2008closed} for details. Indeed, Theorem \ref{thm:Fisher} essentially indicates that
	\begin{equation}\label{GARCH-Fisher}
	I(\theta_0) =  E_{\omega_{\theta_0,2}}^{\theta_0}\left[\frac{1}{2\sigma_1^2} \frac{\partial \sigma_1^2}{\partial \theta}\frac{\partial \sigma_1^2}{\partial \theta^t}{\Big|}_{\theta=\theta_0} \right],
	\end{equation}
	where $\omega_{\theta_0,2}$ is the stationary distribution of $\left\{\left(\sigma_n^2, \frac{\partial \sigma_n^2}{\partial \theta} \right), n\geq 0\right\}$.
	
	We can further link \eqref{GARCH-Fisher} to the formula in \cite{ma2008closed}. For illustration, let us only consider the partial derivative with respect to $\beta$. A direct computation shows that
	\begin{equation*}
	\frac{\partial \sigma_n^2}{\partial \beta} = \sigma_{n-1}^2 + \beta \frac{\partial \sigma_{n-1}^2}{\partial \beta} = \cdots =  \sum_{k=1}^n \beta^{k-1} \sigma_{n-k}^2 + \beta^n \frac{\partial \sigma_0^2}{\partial \beta}.
	\end{equation*}
	(See also equation (7) in \cite{fiorentini1996analytic}.) In addition, under $\omega_{\theta_0,2}$, using the classical procedure of extending the Markov chain to a doubly infinite stationary sequence (see, for example, Section 4 of \cite{bickel1998asymptotic}), we have
	\begin{eqnarray}\label{GARCH-Fisher-rep1}
	\frac{1}{2\sigma_n^2}\left(\frac{\partial \sigma_n^2}{\partial \beta} \right)^2  = \frac{1}{2\sigma_n^2}\left( \sum_{k=1}^\infty \beta^{k-1} \sigma_{n-k}^2 \right)^2,
	\end{eqnarray}
 which has the same distribution as $\frac{1}{2\sigma_0^2} \left( \sum_{k=1}^\infty \beta^{k-1} \sigma_{-k}^2 \right)^2$ as $|\beta|<1$. Combining with
	\eqref{GARCH-Fisher-0} and \eqref{GARCH-Fisher}, we have
	\begin{align}\label{GARCH-Fisher-rep}
	E_{\omega_{\theta_0,2}}^{\theta_0}\left[\frac{1}{\sigma_1^2} \frac{\partial \sigma_1^2}{\partial \beta}\frac{\partial \sigma_1^2}{\partial \beta^t} \right]
		= \lim_{n \rightarrow \infty} E^{\theta_0}\left[\frac{1}{\sigma_n^2} \frac{\partial \sigma_n^2}{\partial \beta}\frac{\partial \sigma_n^2}{\partial \beta^t} \right]
	= E^{\theta_0}\left[ \frac{\left( \sum_{k=1}^\infty \beta^{k-1}\sigma_{-k}^2 \right)^2}{2\sigma_0^2}\right],
	\end{align}
	which is consistent to the closed-form expression provided in (17) of \cite{ma2008closed}.
	
	Similar approach works for the KL-divergence. Let $\theta_i = (\delta_i, \alpha_i, \beta_i)^t$ for $i=0,1$, and denote $\sigma_{n,i}^2$ be the $\sigma_n^2$ evaluated under $\theta_i$. Since \eqref{GARCH-sum} already writes the log likelihood as an additive functional, Theorem \ref{thm:KL} essentially means that
	\begin{eqnarray}\label{KL-GARCH}
	K(\theta_1, \theta_0) =  -\frac{1}{2}E_{\omega_{\theta_1,0}}^{\theta_0}\left[ \log \sigma_{1,1}^2 + \frac{Y_1^2}{\sigma_{1,1}^2} \right]
	 + \frac{1}{2} E_{\omega_{\theta_0,0}}^{\theta_0}\left[
	\log \sigma_{1,0}^2 + \frac{Y_1^2}{\sigma_{1,0}^2}
	\right].
	\end{eqnarray}
	To further link \eqref{KL-GARCH} with the doubly infinite stationary sequence, note that for $i=0,1$, a direct computation leads to
	\begin{align*}
	\sigma_{n,i}^2 & =  \delta_i + \alpha_i Y_{n-1}^2 + \beta_i \sigma_{n-1,i}^2 \\
	& = \delta_i + \alpha_i Y_{n-1}^2 + \beta_i ( \delta_i + \alpha_i Y_{n-2}^2 + \beta_i \sigma_{n-2,i}^2) \\
	& = \cdots \\
	& = \delta_i \sum_{k=0}^\infty \beta_i^k + \alpha_i \sum_{k=0}^\infty \beta_i^k Y_{n-1-k}^2 \\
	& = \frac{\delta_i}{1-\beta_i} + \alpha_i \sum_{k=0}^\infty \beta_i^k Y_{n-1-k}^2,
	\end{align*}
	and so, by stationarity,
	\begin{align*}
	& E_{\omega_{\theta_1,i}}^{\theta_0}\left[ \log \sigma_{1,i}^2 + \frac{Y_1^2}{\sigma_{1,i}^2} \right] \\
	= & \lim_{n \rightarrow \infty}
	E^{\theta_0}\left[ \log \left( \frac{\delta_i}{1-\beta_i} + \alpha_i \sum_{k=0}^\infty \beta_i^k Y_{n-1-k}^2 \right)  + \frac{Y_n^2}{\frac{\delta_i}{1-\beta_i} + \alpha_i \sum_{k=0}^\infty \beta_i^k Y_{n-1-k}^2} \right]\\
	= & E^{\theta_0}\left[ \log \left( \frac{\delta_i}{1-\beta_i} + \alpha_i \sum_{k=0}^\infty \beta_i^k Y_{-(k+1)}^2 \right)  + \frac{Y_0^2}{\frac{\delta_i}{1-\beta_i} + \alpha_i \sum_{k=0}^\infty \beta_i^k Y_{-(k+1)}^2} \right].
	\end{align*}

\begin{remark}
	As demonstrated in \eqref{GARCH-Fisher-rep}, the invariant measure $\omega_{\theta_0,2}$ in Theorem \ref{thm:Fisher} actually incorporates the data of the entire past history. This can be viewed as follows: the invariant probability measure $\omega_{\theta_0,2}$ represents the limiting behaviour of the derivatives of $\ell(\theta; Y_{0:\infty})$, which is equivalent to the limiting behaviour of the derivatives of $\ell(\theta; Y_{-\infty:0})$ when the process is stationary. Similar situation holds for $\omega_{\theta,0}$ in Theorem \ref{thm:KL}.
\end{remark}

To analyze the linear RNN (\ref{6.8a}) and (\ref{RNN5a}), we apply the same method
in \cite{fuh2021a} as follows:
let $W_n=(X_{n-1}, X_{n},Y_{n})^t$ be the Markov chain on ${\cal X}:=({\bf R} \times {\bf R} \times {\bf R})$.	Denote $\eta_n = X_{n-1}^{-1} Y_{n}$ and let
	$\tau_n = (\alpha + \beta \eta_n) \in {\bf R}$.  Let $A_n$ be a $3$-by-$3$ matrix, written  as	
\begin{equation}\label{6.10}
A_{n} = \left( \begin{array}{ccc}
	0 & 1 & 0 \\
	0 & \tau_{n}& 0 \\
	0 & \eta_{n} & 0 \\ \end{array} \right).
\end{equation}
Note that $\{A_n, n \geq 0 \}$ are random matrices driven by the Markov chain $\{W_n, n \geq 0\}$.

Let ${Z}_n=(0,\delta,0)^t \in {\bf R}^{3}$.
Then we have the following state space representation of the linear RNN (\ref{6.8a}) and (\ref{RNN5a}): $W_n$ is a Markov chain govern by
\begin{equation}\label{6.11}
W_{n} = A_{n} W_{n-1} + Z_{n},
\end{equation}
and $Y_n := g(X_n)$, the observed random quantity, is a
non-invertible function of $X_n$.

By Theorem 1 in \cite{fuh2021a}, a sufficient condition for stability is $\alpha + \lambda \beta < 1$, where $\lambda= E_{\Pi} \frac{\mu_{y,1}}{h_{1}}$, with
$\Pi$ as the stationary distribution of the Markov chain $\{(X_n, A_n \cdots A_1), n \geq 0\}$. Condition C1 (the $\omega$-uniformity) can be found in Section 16.5.1 of \cite{Meyn&Tweedie} using the state space representation. Conditions C2--C6 in Theorems \ref{thm:Fisher} and \ref{thm:KL}
hold under the normality assumption in (\ref{6.8a}). By using a similar method as that in
(\ref{LSSMF}), we have a representation of the Fisher information matrix.

In general, a RNN can take as input a variable-length sequence $y=(y_1, \cdots,y_n)$ by recursively processing each symbol while maintaining its internal hidden state $h$. At each time step $n$, the RNN reads the symbol $Y_n \in {\bf R}^q$ and updates its hidden state $h_n \in {\bf R}^p$ by
\begin{eqnarray}\label{RNN1}
h_n = f_{\theta} (Y_n,h_{n-1}),
\end{eqnarray}
where $f_\theta$ is a deterministic non-linear transition function, and $\theta$ is the parameter of $f_\theta$.

For a given RNN model's sequence, by parameterizing a factorization of the joint sequence probability distribution as a product of conditional probabilities, we have
\begin{eqnarray}\label{RNN2}
&~& P(Y_1, \cdots,Y_n) = \prod_{k=1}^n P(Y_k|Y_1,\cdots,Y_{k-1}), \nonumber \\
&~&  P(Y_n|Y_1,\cdots,Y_{n-1} ) =  g_{\theta}(h_{n-1}),
\end{eqnarray}
where $g_\theta$ is a function that maps the RNN hidden state $h_{t-1}$ to a probability distribution over possible outputs, and $\theta$ is the parameter of $g_\theta$.

As noted in \cite{Pascanu1_2014}, given a set of $N$ training sequences $\{y_1^{(n)},\cdots,
y_{T_n}^{(n)}\}$,
the parameters in RNN can be estimated by minimizing the following cost function,
\begin{eqnarray}
J(\theta) = \frac{1}{N} \sum_{n=1}^N \sum_{t=1}^{T_n} d(y_t^{(n)}, g_\theta(h_{t-1}^{(n)})),
\end{eqnarray}
where $d(a,b)$ is a predefined divergence measure between $a$ and $b$, such as Euclidean distance or KL-divergence (or cross entropy).  Theorem \ref{thm:KL}
indicates that the KL-divergence is well-defined in terms of the stationary distribution of
the enlarged Markov chain. This provides a theoretical foundation of using KL-divergence as a cost function in RNN. For the regularization issue, one possible method is the celebrated AIC model selection method in (\ref{eqn:criterion1}) and (\ref{eqn:criterion2}), in which we present a theoretical justification of using this method in RNN.
\end{example}

\begin{example}\label{Ex3}
Temporal Restricted Boltzmann Machine.

A Boltzmann machine is a network with stochastic binary units, which contains a set of visible
units $y \in \{0,1\}^D$ and a set of hidden units $h \in \{0,1\}^P$. The energy of $\{y,h\}$ is defined as
\begin{eqnarray}
E(y,h;\theta) = - \frac{1}{2} y^t L y -  \frac{1}{2} h^t J h  - \frac{1}{2} y^t W h,
\end{eqnarray}
where $\theta =\{W, L, J\}$ are the parameters: $W, L, J$ represent visible-to-hidden,
visible-to-visible, and hidden-to-hidden symmetric interaction terms. The diagonal elements of $L$ and $J$ are set to $0$. The probability that the model assigns to a visible vector $y$ is:
\begin{align}
p(y;\theta) =& \frac{p^*(y;\theta)}{Z(\theta)} = \frac{1}{Z(\theta)}
\sum_h \exp(- E(y,h;\theta)), \\
Z(\theta) =& \sum_y \sum_h \exp(- E(y,h;\theta)),
\end{align}
where $p^*$ denotes unnormalized probability, and $Z(\theta)$ is the partition function (normalizing constant).

Setting both $J=0$ and $L=0$ recovers the well-known restricted Boltzmann machine (RBM) model,
cf. \cite{Hinton_2006}.
 An RBM defines a probability distribution over pairs of vectors, $Y \in \{0,1\}^D$ and
 $H \in \{0,1\}^P$ by the equation
\begin{align}\label{rbm}
p(y,h;\theta) =& P(Y=y, H=h; \theta)  = \frac{1}{Z(\theta)}
 \exp(y^t b_Y + h^t b_H + y^t W h), \\
Z(\theta) =& \sum_y \sum_h \exp(y^t b_Y + h^t b_H + y^t W h),
\end{align}
where $b_Y$ is a vector of bias for the visible vector, $b_H$ is a vector of bias for
the hidden vector, $W$ is the matrix of connection weights, and $Z(\theta)$ is the partition function (normalizing constant).

Next, we consider the temporal restricted Boltzmann machine (TRBM), cf. \cite{Sutskever}. In its simplest form, the TRBM can be viewed as a hidden Markov model (HMM) with an exponentially large state space that
has an extremely compact parameterization of the transition and the emission probabilities.
Denote $a_{1:n}=(a_1,\cdots,a_n)$.
The TRBM defines a probability distribution $P(Y_{1:n}=y_{1:n}, H_{0:n}=h_{0:n})= P((Y_1,\cdots,Y_n)=(y_1,\cdots,y_n), (H_0,\cdots,H_n)=(h_0,\cdots,h_n))$ by the equation
\begin{eqnarray}\label{54}
P(y_{1:n}, h_{0:n}) = \prod_{t=1}^n P(y_t,h_t|h_{t-1}) \bar{\nu} (h_0),
\end{eqnarray}
which has the form as the probability defined in (\ref{eqn:lik}). Here $\bar{\nu}$ can be any suitable
initial distribution of $H_0$. The conditional distribution
$P(Y_t,H_t| h_{t-1})$ is that of an RBM, whose biases for $H_t$ are a function of $h_{t-1}$.
That is
\begin{eqnarray}\label{55}
P(y_t,h_t| h_{t-1}) = \exp( y_t^t b_Y + v^t_t W h_t + h_t^t(b_H + W' h_{t-1})) / Z(h_{t-1}),
\end{eqnarray}
where $b_Y, b_H$ and $W$ are as in Equation (\ref{rbm}), and $W'$ is the weight matrix of the
connection from $H_{t-1}$ to $H_t$, making $b_H + W' h _{t-1}$ be the bias of RBM at time $t$.

Now we need to check conditions C1–C6 in Theorems \ref{thm:Fisher} and \ref{thm:KL} hold under model assumptions in (\ref{54}) and (\ref{55}). First, we note that the state space of $h_t$ is $\mathcal{X} = \{0,1\}^P$, which is finite (although exponentially large); this implies the uniform ergodicity of the underlying Markov chain, and therefore leads that C1 holds. As for the other conditions, note that since the state space $\mathcal{X}$ is finite, all the supremum or integration over $\mathcal{X}$ in C2-C6 are finite. Furthermore, the state space of $y_t$ is $ \{0,1\}^D$, which is finite; this implies that the moment generating function
of $Y$ exists. In addition, since the log of logistic function is infinite differentiable in any local neighbourhood of $\theta$, the supremum over $N_\delta(\theta_0)$ in C2-C6 is finite. This leads that C2-C6 hold.

As noted in \cite{Salakhutdinov_Hinton_2009}, variational learning has the nice property that in addition to trying to maximize the
log likelihood of the training data, it tries to find parameters that minimize the
KL-divergences between the approximating and true posteriors. Theorem \ref{thm:KL}
indicates that the KL-divergence is well-defined in terms of the stationary distribution of
the enlarged Markov chain. This provides a theoretical justification of using the KL-divergence as a cost function in stochastic gradient descent of TRBM.
For the regularization issue, one possible method is the celebrated AIC model selection method in (\ref{eqn:criterion1}) and (\ref{eqn:criterion2}), in which we present a theoretical justification of using this method in TRBM.

However, calculation  of the KL-divergence as well as Fisher information in TRBM and RNN are not straightforward.
This is due to that, unlike the i.i.d. case, the limits in (2) and (3) for TRBM have no explicit form. Traditionally, it is numerically computed by simulating long string of TRBM to approach the limits. Here, thanks to Theorems 1 and 2, we can use Monte Carlo method and other computational technique to evaluate the expectations in (18) and (20) instead. It is worth mentioning that evaluating these expectations are still not straightforward, Theorems \ref{thm:Fisher} and \ref{thm:KL} provide a possible tool for numerical computation via various statistical computation technique.

\end{example}

\section{Conclusion}\label{sec:conclude}
In this paper, we present explicit characterizations of the KL-divergence and
Fisher information for GHMMs, and derive the relationship between these two
important quantities. The results are based on a representation of the log
likelihood and its derivatives  as an additive functional of a MIFS, which
allows one to study the behavior of the log likelihood using SLLN and other
related results. By using these results, we also present the Cram\'{e}r-Rao lower
bound and H\'{a}jek-Le Cam local asymptotic minimax theorem under GHMMs. The
characterization further shows that the KL-divergence in HMM is not convex in
general, which is different from the traditional i.i.d.\ or Markov chain
scenario. Moreover, we provide a theoretical justification of using AIC model selection
in GHMM with finite state space.

It is expected that this representation device will be beneficial
for further studies in \mbox{GHMMs} such as model selection and the generalized method of
moments in stochastic volatility models, exponential tilting estimators and quasi-maximum
likelihood estimators in GHMMs, regularization in RNN with KL-divergence (relative entropy) as the penalized term and other related topics.

\bibliography{references}

\begin{appendices}
\section{Appendix: Proofs of the Main Results}\label{sec:repDeriv}

As Lemmas \ref{thm:rep} and \ref{thm:additiv} are almost the same as Lemmas 3 and 5, respectively, of \cite{fuh_pang_2022} for a two-layer HMM, here we give proofs of these two lemmas in the supplementary for completeness. By using these two lemmas, we first prove Theorems \ref{thm:Fisher}~and~\ref{thm:KL} in Sections \ref{sec:Fisher-HMM}~and~\ref{sec:KL-convergence}, respectively. Then we prove Theorem~\ref{cor:CRLB} and
Corollary~\ref{cor:KL} in Section~\ref{subsec:pf-appl.} based on Theorems \ref{thm:Fisher}~and~\ref{thm:KL}.

\subsection{Proof of Theorem \ref{thm:Fisher}}\label{sec:Fisher-HMM}

To prove Theorem \ref{thm:Fisher}, we only need to prove the following lemma:
\begin{lemma}\label{lemma:Fisher}
Assume conditions C1--C6 hold with $r = 2$. For any $1 \leq j, k \leq q$, we have
\begin{equation*}
    I_{jk}(\theta_0) = -E_{\omega_{\theta_0,2}}^{\theta_0}[g^{\nu(j,k)}(W_{1, \theta_0}^{(2)}, W_{0, \theta_0}^{(2)})],
\end{equation*}
where $E_{\omega_{\theta_0,2}}^{\theta}$ and $g^{\nu(j,k)}$ are the same as in Theorem \ref{thm:Fisher}.
\end{lemma}
Theorem \ref{thm:Fisher} is then a direct consequence of Lemma \ref{lemma:Fisher}.

\begin{proof}
For any $1 \leq j, k \leq q$, by \eqref{Fisher-Bickel}, we have
\begin{equation}\label{pfFisher-1}
    I_{jk}(\theta_0) = -\lim_{n \rightarrow \infty} \frac{1}{n} D^{\nu(j,k)} \ell(\theta_0; Y_{0:n}) ~~~P^{\theta_0} \mbox{-a.s.}
\end{equation}
By \eqref{DL-to-g}, Lemma~\ref{thm:rep}, and the SLLN for Markov random walks in \cite{Meyn&Tweedie}, we have
\begin{align}\label{pfFisher-2}
	& \frac{1}{n}D^{\nu(j,k)} \log L(\theta_0;Y_{0:n}) =  \frac{1}{n} \left\{ \sum_{t=1}^n g^{\nu(j,k)} (W_{t,\theta_0}^{(2)}, W_{t-1, \theta_0}^{(2)}) + g_0^{\nu(j,k)}(W_{0, \theta_0}^{(2)}) \right\}\notag\\
	\rightarrow &~ E_{\omega_{\theta_0,2}}^{\theta_0}[g^{\nu(j,k)}(W_{1, \theta_0}^{(2)}, W_{0, \theta_0}^{(2)})] ~~~P^{\theta_0} \mbox{-a.s.}
\end{align}
Lemma~\ref{lemma:Fisher} is therefore a direct consequence of
\eqref{pfFisher-1} and \eqref{pfFisher-2}.
\end{proof}

\begin{remark}\label{remark:G-to-D2L}
One can further link the function $G$ with the second order derivatives of $\log L(\theta_0; Y_{0:1})$. Define
\begin{equation}\label{G_0}
    G_0(w) = \begin{pmatrix}
    g_0^{\nu(1,1)}(w) &  \cdots & g_0^{\nu(1,q)}(w) \\
    \vdots & \ddots & \vdots \\
    g_0^{\nu(q,1)}(w) & \cdots & g_0^{\nu(q,q)}(w)
    \end{pmatrix}.
\end{equation}
Then, by \eqref{DL-to-g}, we have $$G(W_{1,\theta_0}^{(2)}, W_{0,\theta_0}^{(2)}) = D_\theta^2 \log L(\theta_0; Y_{0:1}) - G_0( W_{0,\theta_0}^{(2)}),$$
and so we have
\begin{equation*}
I(\theta_0) = -E_{\omega_{\theta_0,2}}^{\theta_0}\left[ D_\theta^2 \log L(\theta_0; Y_{0:1}) -G_0(W_{0,\theta_0}^{(2)})\right].
\end{equation*}
\end{remark}

\subsection{Proof of Theorem \ref{thm:KL}} \label{sec:KL-convergence}

To prove Theorem~\ref{thm:KL}, we extend the definition of $W_{n,\theta}^{(r)}$ to $r
\geq 0$ with $W_{n,\theta}^{(0)} = W_{n,\theta}^0$. Note that Lemma \ref{thm:rep} also holds for $r=0$; see \cite{Fuh2006}.
Thus, we can define $\omega_{\theta,0}$ as the stationary distribution of the induced Markov chain
$\{((X_n, Y_n), W_{n,\theta}^{(0)}), n \geq 0\}$.

\begin{proof}[Proof of Theorem \ref{thm:KL}] For simplicity, we prove the case
with $q=1$ and $\theta_1 > \theta_0$; the general case can be proved similarly.

First, it is easy to check that \eqref{DL-to-g} also holds for $r=0$, then we
have
\begin{align}\label{pfKLexist-1}
    & \frac{1}{n} \left[ \log L(\theta_1; Y_{0:n}) - \log L(\theta_0; Y_{0:n}) \right] \notag\\
= & \frac{1}{n} \left\{ \sum_{t=1}^n g^0(W_{t,\theta_1}^{(0)}, W_{t-1,\theta_1}^{(0)})+g_0^0(W_{0,\theta_1}^{(0)})\right\}  - \frac{1}{n}\left\{ \sum_{t=1}^n g^0(W_{t,\theta_0}^{(0)}, W_{t-1,\theta_0}^{(0)})+g_0^0(W_{0,\theta_0}^{(0)})\right\}.
\end{align}
Taking $n \rightarrow \infty$ on both sides of \eqref{pfKLexist-1}, then by
Lemma~\ref{thm:rep}, \eqref{DL-to-g},  and the SLLN for Markov random walks in \cite{Meyn&Tweedie}, we have
\begin{align*}
     K(\theta_1, \theta_0) =  E_{\omega_{\theta_1,0}}^{\theta_1} \left[ g^0(W_{1,\theta_1}^{(0)}, W_{0,\theta_1}^{(0)}) \right]
    - E_{\omega_{\theta_0,0}}^{\theta_1} \left[ g^0(W_{1,\theta_0}^{(0)}, W_{0,\theta_0}^{(0)}) \right],
\end{align*}
which completes the proof for \eqref{eqn:KLrep}.

As for \eqref{HMMKL-converge}, by Taylor expansion, we have
\begin{align}\label{pfKL-1}
    &\ell(\theta_0; Y_{0:n}) - \ell(\theta_1; Y_{0:n}) \notag\\
    = & D^1 \ell(\theta_1; Y_{0:n})(\theta_0-\theta_1) + \frac{1}{2} D^2 \ell(\theta_1; Y_{0:n}) (\theta_0-\theta_1)^2  + \frac{1}{6} \frac{D^3 \ell(\tilde{\theta}_n; Y_{0:n})}{n} (\theta_1-\theta_0)^3,
\end{align}
where $\tilde{\theta}_n \in (\theta_0, \theta_1).$
Dividing both sides of \eqref{pfKL-1} by $n$, we have
\begin{align}\label{pfKL-2}
  &  \frac{1}{n}\left[\ell(\theta_1; Y_{0:n}) - \ell(\theta_0; Y_{0:n})\right] \notag  \\
= &  \frac{D^1 \ell(\theta_1; Y_{0:n})}{n} (\theta_1-\theta_0)  - \frac{1}{2}  \frac{D^2 \ell(\theta_1; Y_{0:n})}{n} (\theta_1-\theta_0)^2  + \frac{1}{6} \frac{D^3 \ell(\tilde{\theta}_n; Y_{0:n})}{n} (\theta_1-\theta_0)^3.
\end{align}

For the first term on the RHS of \eqref{pfKL-2}, by an argument similar to the
proof of Theorem~\ref{thm:Fisher}, we have
\begin{equation}\label{pfKL-1stTerm}
\lim_{n \rightarrow \infty} \frac{D^1 \ell(\theta_1; Y_{0:n})}{n} = E_{\omega_{\theta_1,1}}^{\theta_1}[g^1(W_{1,\theta_1}^{(1)}, W_{0,\theta_1}^{(1)})]  ~~~P^{\theta_1} \mbox{-a.s.}
\end{equation}
At the meantime, we also have $\frac{1}{n} E_x^{\theta_1}[D^1 \ell(\theta_1; Y_{0:n})] = 0$
due to the fact that $L(\theta_1; \cdot)$ is the likelihood under $\theta_1$. Moreover,
\begin{align}\label{pfKL-1stTerm-2}
 \frac{1}{n} E_x^{\theta_1}[D^1 \ell(\theta_1; Y_{0:n})] \rightarrow E_{\omega_{\theta_1,1}}^{\theta_1}[g^1(W_{1,\theta_1}^{(1)}, W_{0,\theta_1}^{(1)})].
\end{align}
 Combining \eqref{pfKL-1stTerm} and \eqref{pfKL-1stTerm-2}, we prove that the first term goes to zero $P^{\theta_1}$-a.s.

For the second term on the RHS of \eqref{pfKL-2}, by \eqref{Fisher-Bickel}, we
have
\begin{equation}\label{pfKL-2ndTerm}
\lim_{n \rightarrow \infty} \frac{D^2 \ell(\theta_1; Y_{0:n})}{n} = -I(\theta_1) ~~~P^{\theta_1} \mbox{-a.s.}
\end{equation}

For the third term in the RHS of \eqref{pfKL-2}, by an argument similar to the
proof of Theorem~\ref{thm:Fisher}, along with the classical uniform LLN (see
\cite{ferguson2017course}, Chapter~16), there exists a constant $C_3>0$ such
that
\begin{align}\label{pfKL-3ndTerm}
\notag
& \limsup_{n \rightarrow \infty} \Bigg\vert \frac{1}{6} \frac{D^3 \ell(\tilde{\theta}_n; Y_{0:n})}{n} (\theta_1-\theta_0)^3 \Bigg\vert \\
\notag
\leq & \limsup_{n \rightarrow \infty} \sup_{\theta \in [\theta_0, \theta_1]} \Bigg\vert \frac{1}{6} \frac{D^3 \ell(\theta; Y_{0:n})}{n} (\theta_1-\theta_0)^3 \Bigg\vert  \\
\leq & C_3|\theta_1-\theta_0|^3 ~~~P^{\theta_1} \mbox{-a.s.}
\end{align}

Thus, taking $n \rightarrow \infty$ on \eqref{pfKL-2} and applying
\eqref{pfKL-1stTerm}--\eqref{pfKL-3ndTerm}, combined with the definition of
$K(\theta_1, \theta_0)$ in \eqref{KL-HMM}, we have
\begin{eqnarray*}
    K(\theta_1, \theta_0) = \frac{I(\theta_1)}{2}(\theta_1-\theta_0)^2 + O\left(|\theta_1-\theta_0|^3\right) =\frac{I(\theta_1)}{2}(\theta_1-\theta_0)^2 + o\left((\theta_1-\theta_0)^2\right),
\end{eqnarray*}
which gives \eqref{HMMKL-converge} when further noticing that $I(\theta_1)
\rightarrow I(\theta_0)$ as $\theta_1 \rightarrow \theta_0$.
\end{proof}

\begin{remark}\label{remark:KL-to-L}
By \eqref{DL-to-g}, we have $g^0(W_{1,\theta}^{(0)}, W_{0,\theta}^{(0)}) = \log L(\theta; Y_{0:1}) - g_0^0( W_{0,\theta}^{(0)})$, so we have
\begin{align*}
    K(\theta_1, \theta_0)  =  E_{\omega_{\theta_1,0}}^{\theta_1} \left[ \log L(\theta_1; Y_{0:1}) - g_0^0(W_{0,\theta_1}^{(0)}) \right] - E_{\omega_{\theta_0,0}}^{\theta_1} \left[ \log L(\theta_0; Y_{0:1}) - g_0^0(W_{0,\theta_0}^{(0)}) \right],
\end{align*}
which further links the KL-divergence to the log likelihood.
\end{remark}

\subsection{Proofs of Theorem \ref{cor:CRLB} and Corollary \ref{cor:KL}}\label{subsec:pf-appl.}

\begin{proof}[Proof of Theorem \ref{cor:CRLB}] Let us begin with \eqref{eqn:CRLB-classical}. For any fixed $n$, the classical
multivariate Cram\'{e}r-Rao lower bound gives
\begin{equation}\label{pf:CRLB-classical}
     E_x^{\theta_0}\left[ \left( v^t(\hat{\theta}_n(Y_{0:n}) - \theta_0) \right)^2\right] \geq v^t I_n^{-1}(\theta_0) v,
\end{equation}
where $I_n(\theta_0) = -E_x^{\theta_0} \left[ D_\theta^2 \log L(\theta_0; Y_{0:n})
\right]$ is the Fisher information based on $Y_{0:n}$.
By using the same argument as
in the proof of Theorem~\ref{thm:Fisher}, we have
\begin{equation}\label{pf:CRLB-Fisher}
    \frac{1}{n} I_n(\theta_0)\rightarrow I(\theta_0).
\end{equation}
Equation~\eqref{eqn:CRLB-classical} immediately follows from \eqref{pf:CRLB-classical}
and \eqref{pf:CRLB-Fisher}.

Let us now turn to \eqref{eqn:CRLB-general}. Denote $P_n^\theta$ as the
probability distribution of $Y_{0:n}$ under $P^\theta$. Then, by Le Cam's
method with squared error loss (\cite{tsybakov2008introduction}, Chapter 2), we
have
\begin{align}\label{pf:CRLB-LeCam}
\inf_{\hat{\theta}_n}  \max_{\theta \in \{\theta_0, \theta_0+\delta v \} } E_x^{\theta}\left[ \Vert \hat{\theta}_n(Y_{0:n}) - \theta  \Vert^2 \right]\geq \frac{\delta^2 \Vert v \Vert^2}{8}\left[ 1 - 2\Vert P_n^{\theta_0} - P_n^{\theta_0 + \delta v} \Vert_{TV}^2 \right],
\end{align}
where $\Vert \cdot \Vert_{TV}$ denotes the total variation distance. In
addition, for any probability distribution $P$~and~$\tilde{P}$, Pinsker's inequality
(\cite{tsybakov2008introduction}, Lemma~2.5(i)) states that $2\Vert P - \tilde{P}
\Vert_{TV}^2 \leq D_{KL}(P \Vert \tilde{P})$, where $D_{KL}(P \Vert \tilde{P}) = \int \log \left(
\frac{P}{\tilde{P}}\right) dP$ is the KL-divergence between $P$ and $\tilde{P}$. By such, we have
\begin{equation}\label{pf:CRLB-Pinskers}
    2 \Vert P_n^{\theta_0} - P_n^{\theta_0 + \delta v} \Vert_{TV}^2
    \leq E_x^{\theta_0 + \delta v} \left[ \log \frac{L(\theta_0+\delta v; Y_{0:n})}{L(\theta_0; Y_{0:n})}\right].
\end{equation}

By using the same argument as in the proof for \eqref{eqn:KLrep} in
Theorem~\ref{thm:KL}, we have
\begin{align} \label{pf:CRLB-Thm2}
   & \frac{1}{n} E_x^{\theta_1} \left[ \log \frac{L(\theta_1; Y_{0:n})}{L(\theta_0; Y_{0:n})}\right] =  \frac{1}{n} E_x^{\theta_1} \left[ \ell(\theta_1; Y_{0:n})\right] - \frac{1}{n} E_x^{\theta_1} \left[ \ell(\theta_0; Y_{0:n})\right] \notag\\
\rightarrow &  E_{\omega_{\theta_1,0}}^{\theta_1}\left[ g^0(W_{1,\theta_1}^{(0)}, W_{0,\theta_1}^{(0)}) \right]
    -    E_{\omega_{\theta_0,0}}^{\theta_1}\left[ g^0(W_{1,\theta_0}^{(0)}, W_{0,\theta_0}^{(0)}) \right] ~~~P^{\theta_1} \mbox{-a.s.} \notag\\
= & K(\theta_1, \theta_0)
\end{align}
for any $\theta_1 \in N_\delta(\theta_0)$. Moreover, similar to the classical
uniform LLN (see \cite{ferguson2017course}), the convergence in
\eqref{pf:CRLB-Thm2} is uniform over any compact subspace of
$N_\delta(\theta_0)$.

Now, recall that $\delta^2 = (n v^t I(\theta_0) v)^{-1}$. By combining
\eqref{HMMKL-converge}, \eqref{pf:CRLB-Pinskers}, and \eqref{pf:CRLB-Thm2}, we
have
\begin{align}\label{pf:CRLB-converge}
\notag
    & \lim_{n \rightarrow \infty} 2 \Vert P_n^{\theta_0} - P_n^{\theta_0 + \delta v} \Vert_{TV}^2
\leq \lim_{n \rightarrow \infty} n K(\theta_0 + \delta v, \theta_0) \\
= & \lim_{n \rightarrow \infty} n \left\{ (\delta v)^t \frac{I(\theta_0)}{2}(\delta v) + o(\Vert \delta v \Vert^2) \right\} = \frac{1}{2}.
\end{align}
Combining \eqref{pf:CRLB-LeCam} and \eqref{pf:CRLB-converge}, we have
\eqref{eqn:CRLB-general} as desired.
\end{proof}

\begin{proof}[Proof of Corollary \ref{cor:KL}]
	For the non-negativeness, note that \eqref{eqn:KLrep} is obtained through
	SLLN for Markov random walks applied to \eqref{KL-HMM}, by which we also have
	\begin{equation*}
	K(\theta_1, \theta_0) = \lim_{n \rightarrow \infty}\frac{1}{n}\left\{ E^{\theta_1}\left[ \ell(\theta_1; Y_{0:n})\right] - E^{\theta_1}\left[ \ell(\theta_0; Y_{0:n}) \right]\right\}.
	\end{equation*}
	However, by Gibb's inequality, we have $E^{\theta_1}\left[ \ell(\theta_1;
	Y_{0:n})\right] - E^{\theta_1}\left[ \ell(\theta_0; Y_{0:n}) \right]  \geq
	0$, so the non-negativeness follows.

	The additivity, on the other hand, is a direct consequence of \eqref{KL-HMM} and the fact that $\log L(\theta_i; Y_{0:n}) = \log L(\theta_i;
	Y_{0:n}^1) + \log L(\theta_i; Y_{0:n}^2)$ for $i=1, 2$ and all $n$.
\end{proof}

\begin{remark}
	The non-negativeness of KL-divergence has also been provided in Section~1.2
	of \cite{douc2011consistency} using a different method. The additivity has been partially investigated in \cite{gorban2010entropy}.
\end{remark}


\newpage

\begin{center}
\Large{
\underline{Supplementary}\\
Kullback-Leibler Divergence and AIC \\in General Hidden Markov Models}

\normalsize{
Cheng-Der Fuh, Chu-Lan Michael Kao and  Tianxiao Pang}
\end{center}

Before proving Lemma \ref{thm:rep}, we need the following definitions.
Since we will differentiate $M_n = {\bf P}_\theta(Y_n) \circ \cdots \circ {\bf P}_\theta(Y_0)$, we need to investigate how the differential operator $D_i$ interacts with the operator $\circ$. Recall ${\bf M}$ and ${\bf P}_\theta$ defined in the first paragraph of Section \ref{subsec:KeyLemma}. Note that for any two given random functions
${\bf P}_\theta(Y_{t+1})$ and ${\bf P}_\theta(Y_{t})$, and any $h_\theta \in {\bf M}$, by
conditions C1--C6 and the dominated convergence theorem, we have
	\begin{align}\label{M-check}
	& D_i \left\{ {\bf P}_\theta(Y_{t}) h_\theta(x) \right\} \notag\\
	= & D_i \left\{
	\int_{s \in \mathcal{X}} p_\theta(s,x) f(Y_t; \theta|x, Y_{t-1})h_\theta (s)Q(ds) \right\} \notag\\
	= &
	\int_{s \in \mathcal{X}} \bigg\{ f(Y_t; \theta|x, Y_{t-1}) h_\theta(s) D_i p_\theta(s,x) + p_\theta(s,x)  h_\theta(s) D_i f(Y_t; \theta|x, Y_{t-1})  \notag\\
	& + p_\theta(s,x) f(Y_t; \theta|x, Y_{t-1}) D_i h_\theta(s) \bigg\} Q(ds)
	\end{align}
and
	\begin{align*}
	\notag
	& D_i \left\{ {\bf P}_\theta(Y_{t+1}) \circ {\bf P}_\theta(Y_{t}) h_\theta(x) \right\} \\
	\notag
	=& 	D_i \bigg\{ \int_{z \in \mathcal{X}} p_\theta(z,x) f(Y_{t+1}; \theta|x, Y_t) \times \\
	& \left( \int_{s \in \mathcal{X}} p_\theta(s,z) f(Y_t; \theta|z, Y_{t-1})h_\theta(s)Q(ds)
	\right) Q(dz)  \bigg\} \\
	\notag
	= &
	\int_{z \in \mathcal{X}} D_i \left\{ p_\theta(z,x) f(Y_{t+1}; \theta|x, Y_t) \right\}  \times \\
	& \left( \int_{s \in \mathcal{X}} p_\theta(s,z) f(Y_t; \theta|z, Y_{t-1})h_\theta(s)Q(ds)
	\right) Q(dz) \\
	\notag
	& +
	\int_{z \in \mathcal{X}} p_\theta(z,x) f(Y_{t+1}; \theta|x, Y_t) \times \\
	& \left( \int_{s \in \mathcal{X}} D_i \left\{ p_\theta(s,z) f(Y_t; \theta|z, Y_{t-1})  h_\theta(s) \right\} Q(ds)\right) Q(dz) \\
	=& \left\{ D_i {\bf P}_\theta(Y_{t+1}) \right\} \circ {\bf P}_\theta(Y_{t}) h_\theta(x) + {\bf P}_\theta(Y_{t+1}) \circ \left\{ D_i ({\bf P}_\theta(Y_{t}) h_\theta (x) )\right\}.
	\end{align*}
By such, we have $D^\nu
\langle M_n \pi \rangle = \langle D^\nu (M_n \pi) \rangle$. Moreover,
for given $\nu_i$ and $\nu_j$, let $\nu_i + \nu_j$ denote the componentwise addition of the vectors. Then, similar to \eqref{M-check}, we have
\begin{align*}
    & D^\nu \left\{ {\bf P}_\theta(Y_{t}) h_\theta(x) \right\} \\
	= & \sum_{\nu_p+\nu_f+\nu_h = \nu} \bigg\{
	\int_{s \in \mathcal{X}} D^{\nu_p} p_\theta(s,x) \times  D^{\nu_f} f(Y_t; \theta|x, Y_{t-1})  D^{\nu_h} h_\theta (s)Q(ds) \bigg\}.
\end{align*}
Therefore, we have
\begin{align}\label{Fuh2006-Lemma3-Deriv}
    & \bigg\vert D^\nu \left\{ {\bf P}_\theta(Y_{t}) h_\theta(x) \right\} -  D^\nu \left\{ {\bf P}_\theta(Y_{t}) g_\theta(x) \right\} \bigg\vert \notag\\
	= & \bigg\vert \sum_{\nu_p+\nu_f+\nu_h = \nu}
	\int_{s \in \mathcal{X}} D^{\nu_p} p_\theta(s,x) D^{\nu_f} f(Y_t; \theta|x, Y_{t-1}) D^{\nu_h} h_\theta (s)Q(ds) \notag\\
	& -  D^{\nu_p} p_\theta(s,x) D^{\nu_f} f(Y_t; \theta|x, Y_{t-1}) D^{\nu_h} g_\theta (s)Q(ds) \bigg\vert.
\end{align}
Hence, if $D^\nu h_\theta(x) \in {\bf M}$ for all $|\nu| \leq r$, then through an argument similar to that in the proof of Lemma 3 in \cite{Fuh2006} (with the condition C1 within replaced by our condition C6), we have $D^\nu \left\{ {\bf P}_\theta(Y_{t}) h_\theta(x) \right\} \in {\bf M}$  for all $|\nu| \leq r$. In addition, by C2 and C3, we have $D^\nu \pi_\theta(x) \in {\bf M}$ for all $|\nu| \leq r$.

We are now ready to prove Lemma \ref{thm:rep}.

\begin{proof}[Proof of Lemma \ref{thm:rep}]
First, based on the argument above, we have $W_n^{(r)} \in {\bf M}^K
:= \{ v=(m_1, \cdots, m_K)^t: m_k \in {\bf M}, 1 \leq k \leq K \}$. This means that $\{((X_n, Y_n), W_n^{(r)}), n \geq 0\}$ is a stochastic process on $(\mathcal{X} \times {\bf R}^d) \times {\bf M}^K$.

To see that $\{((X_n, Y_n), W_n^{(r)}), n \geq 0\}$ is a MIFS, let us investigate the dynamics of $W_n^{(r)}$. Note that for any $\nu_i$,
\begin{align}\label{iteration}
 W_n^{\nu_i}
=&  D^{\nu_i} \big({\bf P}_\theta(Y_n) \circ \cdots \circ {\bf P}_\theta(Y_1) \circ {\bf P}_\theta(Y_0) \big) \notag\\
=&  \sum_{\substack{1 \leq j \leq k \leq K \\
 \nu_i = \nu_j + \nu_k}} \bigg\{ \frac{(\nu_i)!}{(\nu_j)! (\nu_k)!} D^{\nu_k} {\bf P}_\theta(Y_n) \circ  D^{\nu_j} \bigg( {\bf P}_\theta(Y_{n-1}) \circ \cdots  \circ  {\bf P}_\theta(Y_0) \bigg) \bigg\} \notag\\
=& \sum_{\substack{1 \leq j \leq k \leq K \\ \nu_i = \nu_j + \nu_k}} \frac{(\nu_i)!}{(\nu_j)! (\nu_k)!} \left\{ D^{\nu_k} {\bf P}_\theta(Y_n) \circ W_{n-1}^{\nu_j}  \right\}.
\end{align}
Hence, we can define a $K$-by-$K$ {\it matrix form} $A_n = \left\{ a_n^{ij}: 1
\leq i, j \leq K \right\}$, with each $a_n^{ij} \in {\bf M}$ defined as
\begin{equation}\label{aij}
a_n^{ij} = \left\{
\begin{array}{cc}
\frac{(\nu_i)!}{(\nu_j)!(\nu_k)!} D^{\nu_k} {\bf P}_\theta(Y_n) & \mbox{if } \exists 1 \leq k \leq K \mbox{ such~that } \nu_i = \nu_j + \nu_k,\\
0 & \mbox{otherwise}.
\end{array}
\right.
\end{equation}
In addition, for each $K$-by-$K$ ${\bf M}$-valued matrix form $B = \{b_{ij}: 1
\leq i, j \leq K\}$, and each $K$-dimensional ${\bf M}$-valued vector $V =
(V_1, V_2, \cdots, V_K) \in {\bf M}^K$,  we define
\begin{equation}\label{operator}
    B \circ V := \left(
\begin{array}{c}
\sum_{j=1}^K b_{1j} \circ V_j \\
\sum_{j=1}^K b_{2j} \circ V_j \\
\vdots\\
\sum_{j=1}^K b_{Kj} \circ V_j
\end{array}
\right).
\end{equation}
Then by \eqref{iteration}, we have $W_n^{(r)} =A_n \circ W_{n-1}^{(r)}$, and thus
\begin{eqnarray}\label{wn}
W_n^{(r)} = A_n \circ A_{n-1} \circ \cdots \circ A_1 \circ W_{0}^{(r)},
\end{eqnarray}
where $W_{0}^{(r)} = \{ W_0^\nu: |\nu| \leq r \}$ with $W_0^\nu =D^\nu {\bf P}_\theta(Y_0) .$

More importantly, since $W_n^{(r)} =A_n \circ W_{n-1}^{(r)}$, and by \eqref{aij}, the value of $A_n = \left\{ a_n^{ij}: 1
\leq i, j \leq K \right\}$ is determined solely by $Y_n$, we know that the value of $W_n^{(r)}$ is determined solely by $(Y_n, W_{n-1}^{(r)})$. In addition, since the distribution of $Y_n$ is based on $X_n$ and $Y_{n-1}$, and $\{X_n, n \geq 0\}$ is a Markov chain, $((X_n,Y_n), W_n^{(r)})$ is Markovian, as desired.

Finally, for the ergodicity, through a process similar to the proof of Lemma~3 in \cite{Fuh2006}, it can be shown that for $\theta \in N_\delta(\theta_0)$, the MIFS
$\{((X_n, Y_n), W_n^{(r)}), n \geq 0\}$ satisfies Assumption~K in
\cite{Fuh2006}. Furthermore, Lemma~4 in \cite{Fuh2006} holds for the induced
Markov chain $\{((X_n, Y_n), W_n^{(r)}), n \geq 0\}$ on the state space
$(\mathcal{ X} \times {\bf R}^d) \times {\bf M}^K$, which directly leads to
Lemma~\ref{thm:rep}. The proof is completed.
\end{proof}

\begin{remark}
To illustrate {\rm \eqref{wn}}, let $q=1$, i.e., $\theta$ is one-dimensional.
In this case, $\nu \in {\bf R}$ and we can simply label all $|\nu|\leq r$ by
natural order such that $W_n^{(r)}= (W_n^0, W_n^1, \cdots, W_n^r)^t$, the vector
of the first $r$-th derivatives. Then for any $0 \leq k \leq r$, we have
\begin{align*}
W_n^k &= D^k ({\bf P}_\theta(Y_n) \circ \cdots \circ {\bf P}_\theta(Y_1) \circ {\bf P}_\theta(Y_0)) \\
&= \sum_{0 \leq k_1 \leq k} \bigg\{ \frac{k!}{(k_1)! (k-k_1)!} D^{k_1}
{\bf P}_\theta(Y_n) \circ D^{k-k_1} \bigg({\bf P}_\theta(Y_{n-1}) \circ \cdots \circ {\bf P}_\theta(Y_0) \bigg)\bigg\} \\
&=  \sum_{0 \leq k_1 \leq k} C_{k_1}^k \left\{ D^{k_1} {\bf P}_\theta(Y_n) \circ W_{n-1}^{k-k_1}  \right\},
\end{align*}
where $C_a^b = \frac{b!}{a!(b-a)!}$. Therefore  $W_n^{(r)} = A_n \circ
W_{n-1}^{(r)}$ with
\begin{equation}
A_n =
    \left(
    \begin{array}{cccc}
     {\bf P}_\theta(Y_n) & 0 & \cdots & 0 \\
     C_1^1 D^1 {\bf P}_\theta(Y_n) & {\bf P}_\theta(Y_n) & \cdots & 0 \\
     \vdots & \vdots & \ddots & \vdots \\
     C_r^r D^r {\bf P}_\theta(Y_n) & C_{r-1}^r D^{r-1} {\bf P}_\theta(Y_n) & \cdots & {\bf P}_\theta(Y_n)
    \end{array}
    \right),
    \label{An}
\end{equation}
where $0$ denotes the zero function in ${\bf M}$.
\end{remark}

\begin{remark}
Note that $A_n$  in {\rm \eqref{aij}} and $W_n^{(r)}$ in {\rm \eqref{wn}} are
${\bf M}$-valued, other than the traditional ${\bf R}$-valued matrix and vector, respectively. To illustrate this phenomenon, we consider a $D$-state
HMM with one-dimensional parameter $\theta$ case; then $A_n$ in {\rm
\eqref{An}} is a $K$-by-$K$ matrix form with each element being a $D$-by-$D$
matrix (with $0$ being a $D$-by-$D$ zero matrix). In the same manner, although
the operator defined in {\rm \eqref{operator}} appears to be traditional matrix
multiplication, it is different in that the multiplication within each
component is replaced by $\circ$.
Nevertheless, the essential idea is to introduce a matrix form for $W_n^{(r)}$,
which can be used to show that it  forms an ergodic Markov chain via {\rm
\eqref{wn}}.
\end{remark}

\begin{remark}
The critical innovation in this construction is that one needs to consider {\it all} derivatives with order equals or less than $r$ in order to have a Markovian structure. The reason is that, as shown in \eqref{iteration}, the iterated representation of $W_n^{\nu_i}$ involves {\it all} $W_{n-1}^{\nu_j}$ and $W_{n-1}^{\nu_k}$ with $\nu_i = \nu_j + \nu_k$. This is why one will need to consider $W_n^{(r)}$ instead of $W_n^\nu$.

Note that this makes the approach considerably different from previous studies on HMM such as \cite{bickel1998asymptotic}, who study only the second derivatives of $\ell(\theta; Y_{0:n})$ when investigating the Fisher information matrix. We, on the other hand, study {\it all} derivatives with order being equal to or less than two when doing such investigation.
\end{remark}

It is worth mentioning that the feature of getting a neat form in \eqref{wn} is based on a matrix representation \eqref{aij} for all partial derivatives up to the $r$-th order. This largely helps us to obtain the result in Lemma \ref{thm:rep}.

Before proving Lemma \ref{thm:additiv} for general $\nu$, we present a specific form of the first and second order partial derivatives of the log likelihood function as follows. For $|\nu|= 1$, note that we have
\begin{align*}
    \frac{ \langle D^\nu (M_t \pi )\rangle}{\langle M_t \pi \rangle}
     = \frac{ \langle (D^\nu M_t) \pi + M_t(D^\nu \pi) \rangle}{\langle M_t \pi\rangle}  = \frac{ \langle (W_t^\nu) \pi + W_t^0(D^\nu \pi) \rangle}{\langle W_t^0 \pi\rangle}
\end{align*}
for any $t \geq 0$. Therefore
\begin{align}\label{1stderiv}
\notag
	& D^\nu (\log L( \theta; Y_{0:n} ) ) = D^\nu ( \log \langle M_n \pi \rangle)=
	\frac{ \langle D^\nu (M_n \pi )\rangle}{\langle M_n \pi \rangle}\\
\notag
    = & \sum_{t=1}^n \left\{ \frac{\langle D^{\nu}(M_t \pi) \rangle}{\langle M_t \pi \rangle} - \frac{\langle D^{\nu}(M_{t-1} \pi) \rangle}{\langle M_{t-1} \pi \rangle} \right\} + \frac{ \langle D^\nu (M_0 \pi )\rangle}{\langle M_0 \pi \rangle}
	\\
\notag
	 = & \sum_{t=1}^n \bigg\{ \frac{ \langle (W_t^\nu) \pi + W_t^0(D^\nu \pi) \rangle}{\langle W_t^0 \pi\rangle} -
\frac{ \langle (W_{t-1}^\nu) \pi + W_{t-1}^0(D^\nu \pi) \rangle}{\langle W_{t-1}^0 \pi\rangle}
\bigg\} + \frac{ \langle (W_0^\nu) \pi + W_0^0(D^\nu \pi) \rangle}{\langle W_0^0 \pi\rangle} \\
     =: & \sum_{t=1}^n g^\nu ( W_{t}^{(1)}, W_{t-1}^{(1)}) + g_0^\nu(W_0^{(1)}).
\end{align}
That is, the first order derivative of the log likelihood function can be rewritten
as an additive functional of the Markov chain $\{((X_n, Y_n), W_n^{(1)}), n
\geq 0\}$.

To represent the second order partial derivative of the log likelihood, for
$|\nu|=2$, let us write $\nu=\nu_1+\nu_2$ such that $|\nu_1| = |\nu_2|=1$.
Then, we have
\begin{align}
\notag
	& D^{\nu} \log L(\theta;Y_{0:n}) = D^{\nu_1} \left( D^{\nu_2} \log L( \theta ; Y_{0:n}) \right)
	\\
\notag
= & D^{\nu_1} \frac{\langle W_n^{\nu_2} \pi + W_n^0 (D^{\nu_2} \pi)\rangle}{\langle W_n^0 \pi \rangle} \\
\notag
	= & \frac{\langle  W_n^{\nu} \pi + W_n^{\nu_2} (D^{\nu_1} \pi) + W_n^{\nu_1} (D^{\nu_2}\pi) + W_n^0 (D^{\nu}\pi)  \rangle}{\langle W_n^0 \pi \rangle} \\
\notag
	& - \frac{\langle W_n^{\nu_2} \pi + W_n^0 (D^{\nu_2} \pi)\rangle \times \langle W_n^{\nu_1} \pi + W_n^0 (D^{\nu_1} \pi)\rangle}{\langle W_n^0 \pi \rangle^2} \\
    =: & \sum_{t=1}^n g^\nu(W_{t}^{(2)}, W_{t-1}^{(2)}) +  g_0^\nu(W_0^{(2)}),
	\label{2ndderiv}
\end{align}
where
\begin{align}\label{g-for-2ndDL}
  g^\nu(W_t^{(2)}, W_{t-1}^{(2)})= &
\left\{ \frac{\langle  W_t^{\nu} \pi + W_t^{\nu_2} (D^{\nu_1} \pi) + W_t^{\nu_1} (D^{\nu_2}\pi) + W_t^0 (D^{\nu}\pi)  \rangle}{\langle W_t^0 \pi \rangle} \right.
\notag\\
& \left. -
\frac{\langle  W_{t-1}^{\nu} \pi + W_{t-1}^{\nu_2} (D^{\nu_1} \pi) + W_{t-1}^{\nu_1} (D^{\nu_2}\pi) + W_{t-1}^0 (D^{\nu}\pi)  \rangle}{\langle W_{t-1}^0 \pi \rangle}
\right\}
\notag\\
& -
\left\{ \frac{\langle W_t^{\nu_2} \pi + W_t^0 (D^{\nu_2} \pi)\rangle \times \langle W_t^{\nu_1} \pi + W_t^0 (D^{\nu_1} \pi)\rangle}{\langle W_t^0 \pi \rangle^2} \right.\notag\\
& \left. - \frac{\langle W_{t-1}^{\nu_2} \pi + W_{t-1}^0 (D^{\nu_2} \pi)\rangle \times \langle W_{t-1}^{\nu_1} \pi + W_{t-1}^0 (D^{\nu_1} \pi) \rangle}{\langle W_{t-1}^0 \pi \rangle^2} \right\},
\end{align}
and
\begin{align*}
      g_0^\nu(W_0^{(2)})
	= & \frac{\langle  W_0^{\nu} \pi + W_0^{\nu_2} (D^{\nu_1} \pi) + W_0^{\nu_1} (D^{\nu_2}\pi) + W_0^0 (D^{\nu}\pi)  \rangle}{\langle W_0^0 \pi \rangle} \\
	& - \frac{\langle W_0^{\nu_2} \pi + W_0^0 (D^{\nu_2} \pi)\rangle \times \langle W_0^{\nu_1} \pi + W_0^0 (D^{\nu_1} \pi)\rangle}{\langle W_0^0 \pi \rangle^2}.
\end{align*}
That is, the second order derivative of the log likelihood function can also be rewritten
as an additive functional of the Markov chain $\{((X_n, Y_n), W_n^{(2)}), n
\geq 0\}$.

\begin{proof}[Proof of Lemma \ref{thm:additiv}]
We proved the lemma by mathematical induction as follows. As stated in \eqref{1stderiv}, such
$g^\nu$ and $g_0^\nu$ exist for all $|\nu|=1$. Now suppose $g^\nu$ and $g_0^\nu$ exist for all $|\nu|<r$. Then, when $|\nu|=r$, take $\nu_1$ and $\nu_2$ such that $|\nu_2|=1$ and
$\nu_1+\nu_2=\nu$. By induction assumption, $g^{\nu_1}$ and $g_0^{\nu_1}$ exist, therefore we have
\begin{align*}
	 D^{\nu_1+\nu_2} \ell(\theta;Y_{0:n}) = & D^{\nu_2} \left\{ \sum_{t=1}^n g^{\nu_1} (W_t^{(|\nu_1|)}, W_{t-1}^{(|\nu_1|)})  + g_0^{\nu_1}(W_0^{(|\nu_1|)})\right\} \\
	= & \sum_{t=1}^n D^{\nu_2} g^{\nu_1} (W_t^{(|\nu_1|)}, W_{t-1}^{(|\nu_1|)}) + D^{\nu_2} g_0^{\nu_1}(W_0^{(|\nu_1|)}).
\end{align*}
Moreover, as shown in \eqref{1stderiv}, $g^{\nu_1} (W_t^{(|\nu_1|)},
W_{t-1}^{(|\nu_1|)})$ and $g_0^{\nu_1}(W_0^{(|\nu_1|)})$ involve the
derivatives only up to the order of $|\nu_1|$, so $D^{\nu_2} g^{\nu_1}
(W_t^{(|\nu_1|)}, W_{t-1}^{(|\nu_1|)})$ and
$D^{\nu_2}g_0^{\nu_1}(W_0^{(|\nu_1|)})$ involve the derivatives only up to the
order of $|\nu_1|+|\nu_2|=|\nu|$. In other words, they are functions of
$\{W_n^{(|\nu|)}, n \geq 0\}$ as they consist of all derivatives up to the
order of $|\nu|$. Thus, such $g^{\nu}$ and $g_0^{\nu}$ exist for all $|\nu| \leq r$, which completes the proof.
\end{proof}

\begin{remark}
To prove \eqref{DL-to-g} through mathematical induction, we actually only need the exact form of $g^\nu$ and $g_0^\nu$ for $|\nu|=1$ as in \eqref{1stderiv}. However, since the characterization for the Fisher information involves the representation of the second order derivatives in particular, we present the exact form of $g^\nu$ and $g_0^\nu$ for $|\nu|=2$ in \eqref{2ndderiv}.
\end{remark}

\begin{remark}
The representation of $D^\nu \ell(\theta; Y_{0:n})$ also fills the gap in
\cite{Fuh2006}; namely, Section~2.2 in \cite{Jensen2010}, which raises the question of how to deal with the score function and others.
\end{remark}
\end{appendices}

\end{document}